\newcommand{\lV}{\vspace{-1mm}}
\newcommand{\llV}{\vspace{-2mm}}
\newcommand{\lllV}{\vspace{-5mm}}
\newcommand{\mV}{\vspace{1mm}}
\newcommand{\mmV}{\vspace{2mm}}

\newcommand{\xlH}{\hspace{-.5mm}}

\newcommand{\mH}{\hspace{1mm}}
\newcommand{\mmH}{\hspace{2mm}}
\newcommand{\mmmH}{\hspace{5mm}}

\newcommand{\z}[1]{\mathbf{#1}}
\newcommand{\vv}{\mathbf{v}}

\newcommand{\wv}{\mathbf{w}}
\newcommand{\wTv}{\mathbf{w}^T}
\newcommand{\ad}[1]{\text{adj}(#1)}
\newcommand{\av}{\bm{a}}
\newcommand{\ev}{\mathbf{e}}
\newcommand{\xv}{\mathbf{x}}
\newcommand{\yv}{\mathbf{y}}
\newcommand{\bv}{\mathbf{b}}
\newcommand{\zv}{\mathbf{z}}

\newcommand{\N}{\phantom{0}}
\newcommand{\T}[1]{\text{#1}} 

\newcommand{\blu}{\color{blue}}

\newcommand{\mi}{\T{--}} 
\newcommand{\ma}{\T{+}} 

\newcommand{\bproof}{\small{PROOF}. }
\newcommand{\eproof}{$\;_\square$ \normalsize}

\newcommand{\cRan}[3]{#1\leq#2\leq#3}
\newcommand{\iterM}[1]{A^{(#1)}}
\newcommand{\sle}{A\xv=\bv}
\newcommand{\scsle}{A\z{x'}=\z{b'}}
\newcommand{\vecS}[3]{\z{#1}\in{\mathbb{#2}}^{#3}}

\newcommand{\ac}[3]{a^{(#1)}_{#2,#3}}
\newcommand{\hac}[3]{\hat a^{(#1)}_{#2,#3}}

\newcommand{\hA}{\hat{A}}

\newcommand{\yvc}[1]{\yv^{(#1)}}
\newcommand{\yc}[2]{y^{(#1)}_{#2}}
\newcommand{\zvc}[1]{\zv^{(#1)}}
\newcommand{\zc}[2]{z^{(#1)}_{#2}}

\newcommand{\rc}[1]{\rho^{(#1)}}

\newcommand{\w}{\beta_{\max}}
\newcommand{\wb}{\hat{\beta}_{\max}}

\newcommand{\R}{{\mathbb{R}}}
\newcommand{\Q}{{\mathbb{Q}}}
\newcommand{\Z}{{\mathbb{Z}}}

\documentclass[ijoc,nonblindrev]{informs3} %

\OneAndAHalfSpacedXII 

\usepackage{natbib}
 \bibpunct[, ]{(}{)}{,}{a}{}{,}%
 \def\bibfont{\small}%

\usepackage[ruled,vlined]{algorithm2e}

\usepackage{amsmath}
\usepackage{amsfonts}
\usepackage{array}
\usepackage{bm}
\usepackage{boldline}
\usepackage{multirow}
\usepackage{multicol}
\usepackage{float}
\usepackage{enumitem}
\usepackage{color}
\usepackage{tikz}
\usetikzlibrary{matrix,shapes,arrows,positioning,chains}

\usepackage{subfig}

\TheoremsNumberedThrough     

\EquationsNumberedThrough    

\begin{document}

\RUNAUTHOR{A. R. Escobedo} \RUNTITLE{\LaTeX\mH  Starter Document}
\ARTICLEAUTHORS{%
\AUTHOR{Adolfo R. Escobedo}
\AFF{School of Computing and Augmented Intelligence, Tempe, AZ 85281\\ \EMAIL{adres@asu.edu}}
} %

\TITLE{Exact Factorization Updates \\ for Nonlinear Programming}

\RUNTITLE{Exact Matrix Factorization Updates for Nonlinear Programming}

\ABSTRACT{LU and Cholesky matrix factorization algorithms are core subroutines used to solve systems of linear equations (SLEs) encountered while solving an optimization problem. Standard factorization algorithms are highly efficient but remain susceptible to the accumulation of roundoff errors, which can lead solvers to return feasibility and optimality claims that are actually invalid. This paper introduces a novel approach for solving sequences of closely related SLEs encountered in nonlinear programming efficiently and without roundoff errors. Specifically, it introduces rank-one update algorithms for the roundoff-error-free (REF) factorization framework, a toolset built on integer-preserving arithmetic that has led to the development and implementation of fail-proof SLE solution subroutines for linear programming. The formal guarantees of the proposed algorithms are  established through the derivation of theoretical insights. Their advantages are supported with computational experiments, which demonstrate upwards of 75x-improvements over exact factorization run-times on fully dense matrices with over one million entries. A significant advantage of the methodology is that the length of any coefficient calculated via the proposed algorithms is bounded polynomially in the size of the inputs without having to resort to greatest common divisor operations, which are required by and thereby hinder an efficient implementation of exact rational arithmetic approaches.}

 \HISTORY{various sources}

\maketitle\lllV

\section{Introduction}
LU and Cholesky matrix factorization algorithms are widely used in mathematical programming software to solve systems of linear equations (SLEs). They are long established in linear programming (LP), where the factorizations are utilized in the simplex algorithm to move between basic solutions efficiently through the application of factorization updates. They are similarly prevalent in mixed integer linear programming (MILP), where they are used within the branch-and-bound algorithm to efficiently solve LP relaxations of closely related subproblems (e.g., those of close successor nodes). For analogous reasons, LU and Cholesky matrix factorization algorithms are increasingly utilized in nonlinear programming (NLP), where they are deployed within a number of algorithms to solve sequences of SLEs, which are similar to one another but in ways that are substantially different from the LP/MILP setting. In NLP, the updates of interest are rank-one updates and downdates and, more generally, low-rank modifications of the constraint coefficient matrix. As a primary example, rank-one updates of LU and Cholesky factorizations are used to reflect iterative changes to the Hessian and Jacobian matrices associated with the solution of KKT systems of constrained NLPs \citep{sta07eff}. It is worth remarking that low-rank updates have been also prominently featured in LP interior point methods (e.g., see \cite{meh92def,pan20new}).

Low-rank factorization updates are essential for the efficient deployment of various algorithms including fast model predictive control (MPC) \citep{her15dom,kir11fac}, which is employed in the distillation of hazardous chemical compounds and other critical applications \citep{drg17opt}. Rank-one updates of Cholesky factorizations are widely used in a variety of machine learning applications that incorporate linear models or kernels. They are used, for example, to quickly recompute the solution to least squares support vector machine (LS-SVM) classifier problems \citep{fin01eff}, which can be modeled as linearly constrained quadratic programs \citep{oje08low}. They are also utilized to perform Bayesian inference via the expected propagation method \citep{see07bay}, which can be modeled as linearly constrained bilevel NLPs \citep{min13exp}.

Increased numerical stability is among the chief reasons that LU and Cholesky factorization algorithms are preferred over alternative approaches for solving sequences of closely related SLEs. For instance, the Sherman-Morrison-Woodbury (SWM) formula has been widely used (and continues to be, to a certain extent) to perform low-rank modifications of the inverse matrix efficiently \citep{see08low}. However, the implementation of the SWM formula is particularly vulnerable to compounding roundoff errors and their potential for negative consequences. In fact, the benefits of replacing its use with low-rank factorization updates have been corroborated extensively (e.g., see \cite{fin01eff,oje08low}). In spite of their relatively superior numerical stability, matrix factorization algorithms remain susceptible to the accumulation of roundoff errors inherent in floating point computations, which can affect the behavior of optimization solvers and the validity of their outputs. The types of incongruous outcomes that may occur---often unbeknownst to users---include, but are not limited to, optimal solutions being wrongly eliminated from the feasible domain \citep{bai15hig}, a problem being incorrectly identified as having no feasible solutions \citep{pan15var,pur17bou}, and failing to converge to the optimal solution \citep{cho90fur,fin01eff}.

Although these and other invalid outcomes are admittedly infrequent, their non-negligible plausibility detracts from the implicit trust placed on mathematical programming software. Their potential occurrence is specially concerning when obtaining an invalid claim of feasibility/infeasibility or optimality/suboptimality is exceedingly costly or intolerable \citep{mag17cer,sar11rad}. In optimal design applications, for instance, decisions must be taken based on very few observations. Therefore, inferences obtained via the solution of a sequence of linearly constrained NLPs, for which Cholesky factorization rank-one updates can be applied, must be both valid and highly accurate \citep{see07bay}. Similar concerns are germane to the real-time operation of adaptive systems, which are used to protect radio electronic systems (e.g., radar) by processing spatiotemporal signals on a ``sliding-window' \citep{lek18ada}. It is relevant to add that the validity of NLP solver outcomes is essential to address certain theoretical questions. In computational geometry, for example, exact solutions are desired for various core problems that can be cast as convex quadratic optimization programs: Finding the separating hyperplane that maximizes the distance between two polytopes, calculating the distance between two polytopes, and both finding the smallest annulus and finding the smallest ball enclosing a finite set of points in $d$-dimensional space \citep{gar20eff,gle15exa}.

This paper introduces a novel approach for performing exact factorization updates, for the purpose of constructing fail-proof and efficient validation routines for nonlinear programming and other contexts where it is necessary to solve SLEs that undergo low-rank modifications. The featured approach is a direct solution method with various notable differences from existing factorization update algorithms. First, it is founded on specialized integer-preserving arithmetic subroutines, meaning that all operations performed within the featured algorithms---inclusive of  division---are guaranteed to be exact. These subroutines make efficient use of  unlimited-precision data types. In fact, the length of any coefficient they encounter is bounded polynomially in the size of the inputs without having to resort to greatest common divisor (GCD) operations, which are required by any comparable direct solution method approach built on exact rational arithmetic. Another distinctive feature of the proposed approach is that the updated factorization is not obtained by modifying the existing factorization; instead, it is separately reconstructed one row and one column at a time by leveraging a pair of iterative vectors that fundamentally and uniquely connect the two integer-preserving factorizations. The introduction of this approach leads to four main contributions: (1) Derivation of new theoretical insights that help establish the correctness of the factorization update algorithms; (2) Efficient fail-proof algorithms for updating exact LU and Cholesky factorizations in $O(n^2)$ operations; (3) A study to compare the computational performance of the proposed update algorithms with exact refactorization; (4) Additional theoretical and computational analyses that differentiate the proposed algorithms from integer-preserving factorization updates developed for the simplex algorithm. The article is organized as follows. Section \ref{Sec:Background} provides a background on LU and Cholesky factorizations and on the integer-preserving framework featured herein. Section \ref{Sec:Algebra} derives theoretical insights that serve as the foundation of the exact rank-one algorithms, which are introduced in Section \ref{Sec:Main_Alg}. Section \ref{Sec:Tests} introduces experiments on fully-dense matrices that demonstrate the benefits of the proposed algorithms. Finally, Section \ref{Sec:Conclusion} concludes the article and discusses the limitations and future directions of this work.

\section{Background}\label{Sec:Background}
Define an SLE $A\xv=\bv$, with coefficient matrix $A\in\Q^{m\times n}$ of rank $m\le n$, right-hand side parameter vector $\bv\in\Q^m$, and variable vector $\xv\in\Q^n$. An \textit{LU factorization} of $A$ is defined as a pair of matrices $L$ and $U$, the former lower-triangular and the latter upper-triangular, such that $LU=A$. In general, the LU factorization of $A$ is not unique and its computation requires $O(n^3)$ operations. When $A$ is symmetric, it is possible to compute a (unique) factorization of the form $LL^T$ (a lower-triangular matrix, times its transpose) known as the \textit{Cholesky factorization}. Once the respective factorization has been obtained, the SLE can be solved for any instantiation of $\bv$ through \textit{forward substitution} and \textit{backward substitution}, in  $O(n^2)$ operations (fewer operations may be required when $A$ has a special structure). The first of these algorithm entails solving for $\yv\in\Q^n$ in the triangular system $L\yv = \bv$, and the second algorithm entails solving for $\xv\in\Q^n$ in the triangular system $U\xv = \yv$.

The remainder of this section is organized as follows. Section \ref{Subsec:Fac} provides a basic description of SLE rank-one and low-rank updates and one of their primary uses in NLP, and it reviews existing factorization-based algorithms. Section \ref{Subsec:REF_Framework} describes the integer-preserving framework upon which the exact factorizations featured in this work are constructed. Lastly, Section \ref{Subsec:REF-Up} summarizes two approaches for performing LP-related updates on these factorizations, and it explains their inadequacy for performing rank-one updates.

\subsection{Rank-one and Low-rank Factorization Updates}\label{Subsec:Fac}
In mathematical programming and various other fields, it is often necessary to solve a sequence of closely related SLEs. That is, after solving the system $A\xv=\bv$ via factorization (or another suitable method), one needs to solve the updated system $\hA\xv=\bv'$, where $\hA\in\Q^{m\times n}$ has full row-rank and is obtained from a relatively simple modification of $A$, and $\bv'\in\Q^m$ is an arbitrary right-hand side parameter vector. For certain modifications, it is possible to obtain an LU factorization of $\hA$ (i.e., $\hat L\hat U$) in $O(n^2)$ operations by updating an existing LU factorization of $A$. A \textit{rank-one update} is one such modification given by
\begin{equation}\label{eqn:rank-one_SLE}
  \hA = A + \gamma\vv\wv^T,
\end{equation}
where $\gamma\in\Q^1, \vv\in\Q^m$, and $\wv\in\Q^n$ s.t. $\vv\ne\z{0}$ and $\wv\ne\z{0}$ (i.e., the outer product $\vv\wv^T$ has rank one). In the closely related \textit{rank-one downdate}, the matrix $\gamma\vv\wv^T$ is subtracted from $A$. In either case, a symmetry-preserving update can be performed when $\vv=\wv$ and $A$ is symmetric, which is relevant to the Cholesky factorization. Furthermore, Equation \eqref{eqn:rank-one_SLE} can be generalized into a \textit{low-rank} or \textit{rank-$k$ update} by replacing vectors $\vv, \wv$ with matrices $V\in\Q^{m\times k}$, $W\in\Q^{n\times k}$, respectively, where $k\ge1$ is usually much smaller than $m$.

Rank-one updates are core components of NLP algorithms. As a prominent example, the symmetric-rank-one (SR1) formula is used in Quasi-Newton methods to update an approximation matrix $B$ of the Hessian as
\[\hat B= B+\frac{\left(\z{u}-B\z{s}\right)\left(\z{u}-B\z{s}\right)^T}{\left(\z{u}-B\z{s}\right)^T\z{s}},\]
where $\z{u}$ and $\z{s}$ are $n\times 1$ vectors and $B$ is an $n\times n$ matrix (see \citep{wri99num}). The SR1 formula is obtained from the symmetry-preserving rank-one update $A+\gamma\vv\vv^T$ (where $\vv\in\Q^n$ in this case) by setting $A=B$, $\vv=\z{u}-B\z{s}$, and $\gamma=1/[\left(\z{u}-B\z{s}\right)^T\z{s}]$. Moreover, the two-sided rank-one update formula (TR1) is used to update an approximation matrix $C$ of the Jacobian as
\[\hat C= C+\frac{\left(\z{r}-C\z{s}\right)\left(\bm{\mu}^T-\bm{\sigma}^TC\right)}{\bm{\mu}^T\z{s}-\bm{\sigma}^TC\z{s}},\]
where $\z{r}$ and $\bm{\sigma}$ are $m\times1$ vectors, $\bm{\mu}$ is an $n\times1$ vector, and $C$ is an $m\times n$ matrix, with $n\ge m$ (see \citep{gri02con}). The TR1 formula is obtained from Equation \eqref{eqn:rank-one_SLE} by setting $A=C$, $\vv=\z{r}-C\z{s}$, $\wv=\left(\bm{\mu}^T-\bm{\sigma}^TC\right)^T$, and $\gamma=1/(\bm{\mu}^T\z{s}-\bm{\sigma}^TC\z{s}$). Matrices $B$ and $C$ are the coefficients of linearized KKT systems whose solution provides the next point and Lagrangian multipliers used in the optimization algorithm. Accordingly, these SLEs can be efficiently solved via LU factorization updates \citep{sta07eff}.

Various algorithms have been defined for performing low-rank updates efficiently, beginning with \cite{ben65tri} who introduced a rank-$k$ update algorithm for LDU factorizations---which consist of a lower-triangular, a diagonal, and an upper-triangular matrix factor---that iteratively changes matrices $L$ and $U$ by applying Gaussian elimination-type operations (e.g., row reduction, matrix permutation). \cite{fle74mod} developed a variant of this algorithm for the special case when $k=1$, $U=L^T$, and $D$ is positive definite. \cite{kie88num} introduced a related approach for performing rank-one updates on LU factorizations, which may encounter problems when $A$ is rectangular \citep{sta07eff}. \cite{gil74met} describe various approaches for performing rank-one updates on Cholesky factorizations and one such algorithm for matrices that may not be symmetric positive definite. The best performing algorithms utilize plane rotation methods consisting of Householder transformations or products of Givens matrices. Such methods are also commonly utilized for QR factorizations---a decomposition of $A$ into an orthogonal and an upper trapezoidal matrix \citep{ham08upd}.

More recent works on factorization updates have focused on Cholesky factorization algorithms owing to their superior stability on positive-definite and quasi-definite coefficient matrices \citep{den10mul,gil96sta,hig09cho}. In fact, Bennet's algorithm is numerically stable only when $A$ and the rank-$k$ matrices are symmetric and $D$ is positive definite \citep{gil74met}. Yet, the advantageous numerical properties of Cholesky factorization are guaranteed only under certain technical conditions \citep{gil96sta}, not to mention that numerous critical engineering applications deal with coefficient matrices that are neither quasi-definite nor symmetric---e.g., optimal power flow \citep{oh18mul}. To deal with a broader class of matrices, \cite{sta07eff} introduced three LU update algorithms adapted from those of \cite{ben65tri,fle85sta}, and \cite{kie88num}. These adaptations emphasize numerical stability by allowing different row/column permutations not defined in the original versions. Computational results therein demonstrate a superior performance of LU-based methods over QR-based methods including in their application to solve KKT systems of an NLP test set of \cite{hoc80tes}. However, they also demonstrate that both implementations still deviate from the expected theoretical convergence on quadratic optimization problems.

This work derives fail-proof algorithms for performing rank-one updates on exact matrix factorizations. The algorithms are applicable to LU and Cholesky factorizations and, unlike existing LU factorization algorithms, their efficacy does not depend on pivot strategies \citep{hig11gau} or parameter tuning \citep{sta07eff}. The featured algorithms are also applicable to rank-one downdates, which tend to cause problems when using floating-point arithmetic due to possible cancellations when the rank-one matrix is subtracted from $A$ \citep{hig09cho,ols94ran}. To the best of our knowledge, this is the first work to develop an exact direct solution approach for performing NLP-related factorization updates. By ensuring that the outputs of these core subroutines is correct, the proposed algorithms could help unlock important insights across a large number of engineering and scientific applications. The ensuing subsections introduce the exact arithmetic framework that serves as the foundation of the featured factorization-based approach.

\subsection{The Roundoff-error-free Factorization Framework}\label{Subsec:REF_Framework}
The proposed theory and algorithms build on the roundoff-error-free (REF) factorization framework, a direct solution approach that utilizes exact integer-preserving arithmetic. The REF factorization framework includes subroutines for constructing exact LU and Cholesky factorizations and for solving SLEs exactly via REF forward and backward substitution, both for dense \citep{esc15rou} and sparse matrices \citep{lou19exa}. REF factorization subroutines are significantly faster than their exact rational arithmetic counterparts \citep{esc18sol}. In fact, the sparse variants of the REF framework represent the only known exact factorization algorithms for solving SLEs in time proportional to arithmetic work \citep{lou20eff}. They are now included in MATLAB's SuiteSparse libraries \citep{lou20use}, through which they have been employed in real-world engineering applications---e.g., \cite{dia21imp} used the REF framework to improve the accuracy and computational performance of large-scale models for harbor agitation climate assessment. Exact rational arithmetic factorization algorithms generally cannot solve SLEs in time proportional to arithmetic work because they must constantly carry out GCD operations to prevent exponential growth in the \textit{bit-length} or encoding size of the matrix entries \citep{web19sol}. Conversely, the sparse REF factorization algorithms achieve the former guarantee in part through a set of special properties derived from the integer-preserving Gaussian elimination algorithm (IPGE).

To establish a proper foundation for this work, the ensuing paragraphs introduce assumptions,  notational conventions, and the IPGE and REF factorization algorithms.

\begin{assumption}\label{assu:AMatrix}
\textnormal{For the remainder of this work, let $A\in{\mathbb{Z}^{n\times n}}$ and $\hA = A + \gamma\vv\wv^T$ be nonsingular, and assume that $\bv,\vv,\wv\in\Z^n$  and $\gamma\in\Z^1$}.
\end{assumption}

Let $\iterM{k}$ be the $k$th-iteration matrix of IPGE, for integer $\cRan{0}{k}{n}$, where $\iterM{0}:=A$ and denote its individual entries as $\ac{k}{i}{j}$, for $\cRan{1}{i,j}{n}$. Additionally, let $\rc{k}$ denote the pivot element selected from $\iterM{k\mi1}$ to perform the $k$th iteration of IPGE, where $\rc{0}:=1$.

The recursive formula for calculating IPGE entry $\ac{k}{i}{j}$ is given by:
\begin{eqnarray}
\ac{k}{i}{j} = \left\{
\begin{array}{lll}
\ac{k-1}{i}{j} & \N & \text{ if } i=r_k \\
\left(\rc{k}\ac{k-1}{i}{j}-\ac{k-1}{r_k}{j}\ac{k-1}{i}{c_k}\right)/\rc{k\mi1} & \N & \text{ otherwise}
\end{array}
\right.\qquad \text{ for } \; k=1\ldots n\label{eqn:ipge}
\end{eqnarray}
where $\cRan{1}{k}{n}$, $\cRan{k}{i,j}{n}$, and $\rc{k}\neq0$ for all $k$; and where $r_k$ and $c_k$ are the row and column indices, respectively, of $\rc{k}$ in $A^{k\mi1}$. Note that this represents the Gaussian elimination version of IPGE, which outputs an echelon form (i.e., upper-triangular) matrix and is sufficient to obtain the REF factorizations. The Gauss-Jordan elimination version, which outputs a reduced echelon form (i.e., diagonal) matrix, is obtained by setting $\cRan{1}{i,j}{n}$ in each step and, consequently, entails more computational effort.

\begin{assumption}\label{assu:pivot}
\textnormal{Fix $\rc{k}=\ac{k\mi1}{k}{k}\neq0$ (i.e., $r_k=c_k=k)$, for $k\geq1$.}
\end{assumption}

Prior to introducing the REF factorization algorithms, it is useful to elaborate on the  preceding assumptions. Assumption \ref{assu:AMatrix} does not lead to a loss of generality since any rational matrix can be multiplied by the lowest common denominator of its entries and any  finite-precision matrix can be multiplied (i.e., right-shifted) by an adequate power of 10 to yield an SLE whose coefficients are all integers. Assumption \ref{assu:pivot} implies that IPGE does not need to perform row/column permutations to find nonzero pivot elements and is adopted for simplicity. For extended algorithms that deal with the possibility of encountering zeros along the diagonal, see \cite{esc16fou}.

Next, let $[k]$ be shorthand for the ordered index set $\{1,\dots,k\}$, where $k\ge1$. Expanding from this notation, denote $A^{[k],j}_{[k],i}$ as the submatrix induced by the ordered column-index set $\{1,\dots,k,j\}$ and the ordered row-index set $\{1,\dots,k,i\}$ of $A$; similarly, denote $\vv_{[k],i}$ as the subvector induced from $\vv$ by the ordered index set $\{1,\dots,k,i\}$. This notation is utilized throughout this paper, and it is helpful in this section for presenting the three key properties of IPGE. First, every division performed in the algorithm is guaranteed to be exact \citep{bar68syl}, that is, each dividend is an integer multiple of its divisor. Second, the maximum bit-length of any IPGE entry, denoted as $\w$, is bounded polynomially as $\beta_{\max}\leq \lceil n\log(\sigma\sqrt{n})\rceil$, where $\sigma:=\underset{i,j}{\max\mH}|\ac{0}{i}{j}|$ \citep{bar72com}. When $A$ is symmetric positive definite, the bound on $\beta_{\max}$ tightens to  $\lceil n\log(\sigma)\rceil$; and when $A$ is sparse, it is reduced to $\lceil n\log(\sigma\sqrt{\delta})\rceil$, where $\delta$ is the minimum between the number of non-zeros in the most dense row and the most dense column of $A$ \citep{lou19exa}. Third, each IPGE entry $\ac{k}{i}{j}$ equates to a specific subdeterminant of $A$, as follows \citep{edm67sys}:
\begin{eqnarray}
\ac{k}{i}{j} = \left\{
\begin{array}{lll}
              (-1)^{i+k}\det\left(A_{\{[k]\backslash i\}\cup j}^{[k]}\right) & \N & \text{ if } i\leq k \mmV\\
                \det\left((A)_{[k],j}^{[k],i}\right) & \N & \text{ otherwise};
\end{array}\label{eqn:subdet}
\right.
\end{eqnarray}
for $0\leq{k}\leq{n}$ and $k\leq{i,j}\leq{n}$. Based on this characterization, Assumption \ref{assu:pivot} implies that, for any $k\ge1$, the subvectors $A^{1}_{[k]},A^{2}_{[k]}\dots, A^{k}_{[k]}$ are linearly independent (since $\rc{k}\ne0$ is the $k$th leading principal minor of nonsingular $A$).

The REF LU factorization of $A$, henceforth abbreviated as REF-LU($A$), is an LD$^{\mi1}$U factorization, that is, it consists of a lower-triangular matrix $L$, the inverse of a diagonal matrix $D$, and an upper-triangular matrix $U$. The contents of the three matrix factors are obtained from the iterative entries of the Gaussian
elimination version of the IPGE algorithm applied to $A$ and are given by:
\def\arraystretch{2}\arraycolsep=3pt
\begin{eqnarray}
l_{i,j} &= \ac{j\mi1}{i}{j}, &\T{ for } i\ge j;\label{eqn:reflu1}\\
d_{i,i} &= \rc{i\mi1}\rc{i}=\ac{i\mi2}{i\mi1}{i\mi1}\ac{i\mi1}{i}{i}, & \T{ for  all } i; \T{ and}\label{eqn:reflu2}\\
u_{i,j} &= \ac{i\mi1}{i}{j}, &\T{ for } i\le j;\label{eqn:reflu3}
\end{eqnarray}
where $\cRan{1}{i,j}{n}$. However, $D$ does not need to be stored. Its entries can be generated from the diagonal of $L$ or $U$---specifically, $\ac{i\mi2}{i\mi1}{i\mi1}=l_{i\mi1,i\mi1}=u_{i\mi1,i\mi1}$ and $\ac{i\mi1}{i}{i}=l_{i,i}=u_{i,i}$,
where it is assumed that $\ac{\mi1}{0}{0}:=\rc{0}=1$. When $A$ is symmetric, $U=L^T$, thereby inducing the REF Cholesky factorization, whose exact expression is given by  $(LD^{\mi1/2})(LD^{\mi1/2})^T$.

REF forward substitution on a vector $\bv$ is performed with the lower-triangular factor $L$ of REF-LU($A$) by initializing and then iteratively updating vector $\vecS{y}{Z}{n}$, for iterations $k=0,\dots,n\mi1$, as follows:
\begin{eqnarray}
y_i=\left\{\def\arraystretch{.75}\arraycolsep=3pt
\begin{array}{lll}
b_i                      & \N & \text{ if } k=0,\mV\\
l_{1,1}y_i - l_{i,1}y_1& \N & \text{ if } k=1,\mV\\
\frac{(l_{k,k}y_i - l_{i,k}y_k)}{l_{k\mi1,k\mi1}}& \N & \text{ otherwise }
\end{array}
\right.; \qquad \text{ for } \; i=k\ma1\ldots n.\label{eqn:intfsB}
\end{eqnarray}
The output vector is the solution to the SLE $LD^{\mi1}\yv=\bv$ or, equivalently, $\yv=(LD^{\mi1})^{-1}\bv$. However, to solve $\sle$ in full without roundoff errors, REF forward and backward substitution must be applied on the scaled SLE $\scsle$, where $\xv':=\det(A)\xv$ and $\bv':=\det(A)\bv$. The forward substitution vector $\z{y}'$ for the scaled system $LD^{\mi1}\yv'=\bv'$ is equivalently obtained without roundoff errors by evaluating \eqref{eqn:intfsB} and setting $\yv':=\det(A)\yv=l_{n,n}\yv$.

Having obtained $\yv'$, REF backward substitution is performed as follows:
\begin{flalign}
x'_i = \frac{1}{u_{i,i}}\left(y'_i - \underset{j=i+1}{\overset{n}{\sum}}u_{i,j}x'_j\right)\hspace{.5in}\text{ for } \; i=n\dots1. \label{eqn:intbs}
\end{flalign}
Afterward, the exact solution to the original SLE can reported to any desired precision through the equation:
\begin{eqnarray}
x_i=\frac{x'_i}{\det(A)}=\frac{x'_i}{l_{n,n}} \hspace{.5in}\text{ for } \; i=1\dots n.\label{eqn:REFsol}
\end{eqnarray}

Before proceeding, it is important to state that the worst-case computational complexities (WCC) of IPGE, REF factorization, and REF substitution are as follows:
\begin{flalign}
\text{WCC(IPGE/REF factorization)}&=O(n^3(\w\log \w \log\log\w))\label{eqn:wcc1}\\
\text{WCC(REF substitution)}\hspace{1.15cm}&=O(n^2(\w\log\w\log\log\w)).\label{eqn:wcc3}
\end{flalign}
In Equations \eqref{eqn:wcc1} and \eqref{eqn:wcc3}, the WCC measures use the maximum bit-length $\w$ from the IPGE algorithm to account for the added complexity of operand growth in the exact factorization's entries. The expression in the innermost parentheses of both equations represents the cost of multiplying/dividing two integers of bit-length $\w$ according to FFT techniques \citep{sch71sch, knu81art}; the quantity outside the innermost parentheses represents the algorithms' number of operations.

\subsection{Updating the REF Factorizations}\label{Subsec:REF-Up}
\cite{esc17rou} introduced algorithms for performing various LP-related updates on the REF LU and Cholesky factorizations, namely, addition, deletion, and replacement of a single row/column of $A$. The \textit{push-and-swap} column replacement approach developed therein contrasts with the traditional \textit{delete-insert-reduce} update approach, although they both require $O(n^2)$ operations. The latter approach immediately deletes the exiting column, inserts the incoming column in a strategic position, and performs row-reduction operations to return the factors to triangular form (for a survey of various update algorithms that can be categorized under the delete-insert-reduce approach, we refer the reader to \cite{elb12rev}). Applying this traditional approach on the REF factorizations leads to a loss in the information that was used to guarantee exact divisibility during each iteration of the factorization process. Without this information, the IPGE pivoting process must be restarted to guarantee exact divisibility in the new row reduction operations, causing prohibitive increases in computational effort. In a set of computational experiments on fully dense matrices performed in \cite{esc17rou}, the run-times of this update approach always exceeded those required by the full factorization; in fact, they were over 90x slower for basis matrices with dimension $n=500$. Conversely, the push-and-swap column replacement approach preserves the special structure of the REF factorization by repeatedly permuting the column exiting the basis with its right-adjacent column until it is pushed out of the factorization and swapped with the (updated) incoming column. This special procedure avoids growth in the encoding size of the matrix entries, specifically, it ensures the entries of the updated factorization retain the IPGE bound on $\beta_{\max}$ (see Section \ref{Subsec:REF_Framework}).

A crucial distinction of a rank-one update is that most, if not all, columns of $A$ change at once rather than a single row or column, when the operation defined in \eqref{eqn:rank-one_SLE} is applied to yield $\hA$. This means that, while REF-LU($\hA$) could be obtained as a sequence of column replacements, doing so would require $O(n^3)$ operations---$O(n^2)$ operations for each of $O(n)$ column updates. This would cancel out the operations savings expected of the factorization update and, therefore, a fundamentally different approach is needed to perform  efficient rank-one updates on the REF factorizations. Inefficient algorithms would also result from a direct adaptation of other update algorithms (e.g., \cite{gil87mai,sta07eff}), that is, through the replacement of their floating-point operations with exact arithmetic. In greater detail, because the division operations that would be involved are not guaranteed to be exact, such implementations would entail switching to exact rational arithmetic. This would effectively eliminate the advantages of integer-preserving arithmetic vis-\'a-vis the latter methodology of exact computation, which include faster run times and lower memory requirements. We direct the reader to \cite{esc18sol,lou19exa} for comparisons of these two exact methodologies for LU factorization and forward/backward substitution on dense and sparse matrices.

Conversely, the proposed algorithms can be utilized to perform a column (or row) replacement in $O(n^2)$ operations, since this matrix modification can  be cast as a rank-one update. This additional use is analyzed and evaluated computationally in Section \ref{Sec:Tests}.

\section{Theoretical Insights}\label{Sec:Algebra}
This subsection derives theoretical insights that have special import with the REF LU and Cholesky factorization algorithms. The ensuing discussion focuses on the simpler update $\hA = A + \vv\wTv$ without loss of generality, since the scaled outer product $\gamma\vv'\wTv$ can be expressed as $\vv\wTv$, where $\vv:=\gamma\vv'$ and $\vv'\in\mathbb{Z}^{n}$. The next two theorems introduce new identities of the adjoint matrix, adj$(\cdot)$. An ensuing corollary connects these identities to the IPGE algorithm, from which the entries of the REF-LU factorization are obtained (see Section \ref{Subsec:REF_Framework}). Although this work assumes that $A\in\mathbb{Z}^{n\times n}$, $\vv,\wv\in\mathbb{Z}^{n}$, and $\gamma\in\Z^1$, the theoretical results presented in this section extend to rational- and real-numbered matrices.

\begin{theorem}\label{thm:Adju}
Let $A$ be nonsingular. Then, the following identity holds:
\begin{equation*}
\ad{A}\vv=\ad{A+\vv\wTv}\vv.
\end{equation*}
\end{theorem}
\bproof The adjoint matrices of $A$ and $A+\vv\wTv$ are related through the equation \citep{els81eig}:
\begin{eqnarray}\label{eqn:AdjRel}
&\ad{A+\vv\wTv}=\ad{A}+\wTv\ad{A}\vv A^{-1}-\ad{A}\vv\wTv A^{-1}\\\label{eqn:AdjRel_a}
\Rightarrow&\ad{A+\vv\wTv}-\ad{A}=\left[\wTv\ad{A}\vv\;I_n -\ad{A}\vv\wTv\right]A^{-1},\label{eqn:AdjRel_b}
\end{eqnarray}
where $I_n$ is the identity matrix of order $n$. Multiplying by $\vv$ from the right gives:
\begin{eqnarray*}
\ad{A+\vv\wTv}\vv-\ad{A}\vv=\left[\wTv\ad{A}\vv\;I_n -\ad{A}\vv\wTv\right]A^{-1}\vv.
\end{eqnarray*}
Therefore, the desired result is established by demonstrating that
\begin{eqnarray}
&\left[\wTv\ad{A}\vv\;I_n -\ad{A}\vv\wTv\right]A^{-1}\vv=\z{0}_n\label{eqn:AdjRel2}\\
\Leftrightarrow &\left[\wTv\ad{A}\vv\;I_n -\ad{A}\vv\wTv\right]\ad{A}\vv=\z{0}_n \label{eqn:AdjRel3}
\end{eqnarray}
where $\z{0}_n$ is the zero-vector of size $n$. This is shown by redistributing the left-hand side of \eqref{eqn:AdjRel3} as:
\begin{equation}\label{eqn:AdjRel4}
\left[\wTv\ad{A}\vv\right]\ad{A}\vv -\left[\ad{A}\vv\wTv\right]\ad{A}\vv
= \left[\wTv\ad{A}\vv\right]\ad{A}\vv -\ad{A}\vv\left[\wTv\ad{A}\vv\right]
=\z{0}_n,
\end{equation}
where the first equality in \eqref{eqn:AdjRel4} uses the associativity of matrix multiplication.\eproof

\begin{theorem}\label{thm:AdjuT}
Let $A$ be nonsingular. Then, the following identity holds:
\begin{equation*}
\wTv\ad{A}=\wTv\ad{A+\vv\wTv}.
\end{equation*}
\end{theorem}
\bproof Similar to Theorem \ref{thm:Adju}, we begin with the difference between the adjoint matrix of $A$ and of its rank-one update (see \eqref{eqn:AdjRel_b}), which multiplied from the left by $\wTv$ gives:
\begin{eqnarray*}
\wTv\ad{A+\vv\wTv}-\wTv\ad{A}=\wTv\left[\wTv\ad{A}\vv\;I_n -\ad{A}\vv\wTv\right]A^{-1}
\end{eqnarray*}
In this case, the desired result is established by demonstrating that
\begin{eqnarray}
&\wTv\left[\wTv\ad{A}\vv\;I_n -\ad{A}\vv\wTv\right]A^{-1}=\z{0}_n\\
\Leftrightarrow &\wTv\left[\wTv\ad{A}\vv\;I_n -\ad{A}\vv\wTv\right]\ad{A}=\z{0}_n.
\end{eqnarray}
This is shown by reorganizing the left-hand side of the latter equation, culminating with the expression:
\begin{equation}
\wTv\left[\wTv\ad{A}\vv\right]\ad{A} -\left[\wTv\ad{A}\vv\right]\wTv\ad{A}=\z{0}. \mmH_{\square}
\end{equation}\normalsize

\begin{corollary}
Let $\xv'$ denote the vector obtained from performing REF forward substitution (see \eqref{eqn:intfsB}), followed by REF backward substitution (see \eqref{eqn:intbs}) with REF-LU($A$) on $\vv':=\det(A)\vv$. Additionally, let $\hat{\xv}''$ denote the vector obtained by performing REF forward substitution, followed by REF backward substitution, with REF-LU($\hA$) on $\vv'':=\det(\hA)\vv$. It must be the case that $\xv'=\hat{\xv}''$.
\end{corollary}
\bproof From the given information, $\xv'$ and $\hat{\xv}''$ satisfy the respective SLEs $A\xv'=\vv'$ and $\hA\hat{\xv}''=\vv''$. Based on the properties of REF backward substitution, we have that:
\[\xv'=\det(A)\xv=\det(A)A^{-1}\vv=\ad{A}\vv=\ad{\hA}\vv=\det(\hA)\hA^{-1}\vv=\det(\hA)\hat{\xv}=\hat{\xv}''.\mmH_{\square}\]
The ensuing theorems extend the implications of this result, which in and of itself is insufficient for reconstructing REF-LU($\hA$) from REF-LU($A$). To continue, it is convenient to state a basic identity. \eproof

\begin{proposition}\label{prop:nesting}
For any nonsingular lower-triangular matrix $\Lambda\in\mathbb{R}^{n\times n}$,
\begin{equation}
  \left(\Lambda^{-1}\right)^{[k]}_{[k]}=\left(\Lambda^{[k]}_{[k]}\right)^{-1}.
\end{equation}
That is, the first $k$ rows and columns of $\Lambda^{-1}$ are exactly the inverse of the submatrix induced by the first $k$ rows and columns of $\Lambda$.
\end{proposition}

\begin{theorem}\label{thm:Inverse_Identity}
Let $LD^{-1}U$ and $\hat L\hat D^{-1}\hat U$ be the REF-LU factorizations of $A$ and $\hA=A+\vv\wTv$, respectively. The result of applying forward substitution on $\vv$ using $L$ matches the result of applying forward substitution on $\vv$ using $\hat L$, that is,
\begin{equation*}
(LD^{-1})^{-1}\vv=(\hat L\hat D^{-1})^{-1}\vv.
\end{equation*}
\end{theorem}
\bproof The adjoint matrix $\ad{A}$ can be re-expressed using REF-LU($A$) as:
\begin{flalign*}
\ad{A}&=\det(A)A^{-1}\\
&=\det(A)\left(LD^{-1}U\right)^{-1}\\
&=\det(A)\mH U^{-1}\left(LD^{-1}\right)^{-1}.
\end{flalign*}
From these equations, the $n$th row of $\ad{A}$, written succinctly as $(\ad{A})^{[n]}_n$, is given by
\begin{flalign*}
(\ad{A})^{[n]}_n&=\det(A)\mH (U^{-1})^{[n]}_n\left(LD^{-1}\right)^{-1}\\
&=\det(A)\mH \left[0\mH 0 \dots 0\mmH 1/u_{m,m}\right]\left(LD^{-1}\right)^{-1}\\
&=u_{m,m}\mH \left[0\mH 0 \dots 0\mmH 1/u_{m,m}\right]\left(LD^{-1}\right)^{-1}\\
&=\ev^T_n\left(LD^{-1}\right)^{-1}\\
&=\left(\left(LD^{-1}\right)^{-1}\right)^{[n]}_n;
\end{flalign*}
where $\ev_n$ is the $n$th elementary vector of length $n$. To extend this analysis to characterize an arbitrary row of $\left(LD^{-1}\right)^{-1}$, let $L(A^{[k]}_{[k]})$ and $D^{-1}(A^{[k]}_{[k]})$ be the lower-triangular matrix and diagonal matrix corresponding to REF-LU($A^{[k]}_{[k]}$), that is, the REF LU factorization of the (nonsingular) submatrix induced by the first $k$ rows and columns of $A$. Since the product $LD^{-1}$ is nonsingular and lower-triangular, the inverse matrices $\left(L(A^{[k]}_{[k]})D^{-1}(A^{[k]}_{[k]})\right)^{-1}$ nest atop one another as $k$ increases  per Proposition \ref{prop:nesting}, meaning that the following relationship holds, for $k=1,\dots,n$:
\begin{equation}
  (\left(LD^{-1}\right)^{-1})^{[k]}_{[k]}=\left(L(A^{[k]}_{[k]})D^{-1}(A^{[k]}_{[k]})\right)^{-1}.
\end{equation}
From the preceding series of equations, this implies that the $k$th row of matrix $\left(LD^{-1}\right)^{-1}$ can be written as
\[\left(\left(LD^{-1}\right)^{-1}\right)^{[k]}_k=\left[\left(\ad{A}^{[k]}_{[k]}\right)^{[k]}_k \mH \z{0}^T_{n-k}\right],\]
where $\z{0}^T_{n-k}$ is the (row) 0-vector of dimension $n\mi k$. Lastly, applying Theorem \ref{thm:Adju} gives that:
\[\left(\ad{A}^{[k]}_{[k]}\right)^{[k]}_k\vv=\left(\ad{\hA}^{[k]}_{[k]}\right)^{[k]}_k\vv,\]
for $k=1,\dots,n$, thereby establishing the desired result. \eproof

\begin{theorem}\label{thm:Inverse_Identity_b}
Let $L D^{-1} U$ and $\hat L\hat D^{-1}\hat U$ be the REF-LU factorizations of $A$ and $\hA=A+\vv\wTv$, respectively. The result of applying forward substitution on $\wv$ using $U^T$ matches the result of applying forward substitution on $\wv$ using $\hat U^T$, that is,
\begin{equation*}
(U^TD^{-1})^{-1}\wv=(\hat U^T\hat D^{-1})^{-1}\wv.
\end{equation*}
\end{theorem}
\bproof The REF LU factorizations of $A^T$ and $\hA^T$ are exactly the transpose of the factorizations REF-LU($A$) and REF-LU($\hA$), respectively \citep{esc17rou}. This means that $U^T$ and $\hat U^T$ are the corresponding lower triangular matrices needed to perform REF forward substitution on $\wv$ with REF-LU($A^T$) and REF-LU($\hA^T$). Now, from Theorem \ref{thm:AdjuT}, we have that:
\begin{flalign*}
&\wTv\ad{A}=\wTv\ad{\hA}.\\
\Leftrightarrow\mmH&\ad{A}^T\wv=\ad{\hA}^T\wv.\\
\Leftrightarrow\mmH&\ad{A^T}\wv=\ad{\hA^T}\wv.
\end{flalign*}
Therefore, the desired result is obtained through a parallel line of arguments as Theorem \ref{thm:Inverse_Identity}. \eproof

Theorems \ref{thm:Inverse_Identity} and \ref{thm:Inverse_Identity_b} prove that performing REF forward substitution on $\vv$ (on $\wv$, resp.) with REF-LU($A$) or with REF-LU($\hA$) (with REF-LU($A^T$) or with REF-LU($\hA^T$), resp.) yields identical results. It is necessary to go one step further and show that the intermediary update vectors calculated during this algorithm, with either the original or updated factorizations, also match. To that end, we introduce a stepwise recursion of REF forward substitution in Algorithm \ref{FS_Step}, denoted succinctly as REF-FS-Step. In words, Algorithm \ref{FS_Step} receives the forward substitution vector from step $k\mi1$ (say $\yvc{k\mi1}$), column $k$ of $L$ (say $L^{k}_{[n]}$), the pivot used in step $k\mi1$ ($\rc{k\mi1}=l_{k\mi1,k\mi1}$), and the iteration counter ($k$); it uses these inputs to calculate and return the forward substitution vector for step $k$, where $\cRan{1}{k}{n\mi1}$. Therefore, running the original REF forward substitution algorithm on the update vector $\vv$ with REF-LU($A$) is equivalent to performing the recursion $\yvc{k}=$ REF-FS-Step($\yvc{k\mi1},L^{k}_{[n]},l_{k\mi1,k\mi1},k$), for $k=1,\dots,n\mi1$, where $\yvc{0}=\vv$.

\begin{algorithm}
\SetKwInOut{Input}{input}\SetKwInOut{Output}{output}\SetKw{KwRet}{return}
\Input{$\yvc{k\mi1}$, $L^{k}_{[n]}$, $\rc{k\mi1}$, $k$}

    {\bf let} $\yvc{k}\in\mathbb{Z}^{n}$\\
    \For{$i=k+1,\dots,|\yvc{k\mi1}|$}
    {
       $\yc{k}{i}=l_{k,k} \yc{k\mi1}{i} - l_{i,k}\yc{k\mi1}{k}$

        \If{$k>0$}
        {
            $\yc{k}{i} \leftarrow\yc{k}{i}/\rc{k\mi1}$
        }
    }

\KwRet{$\yvc{k}$}
\caption{REF Forward Substitution Recursive Step (REF-FS-Step)\label{FS_Step}}
\end{algorithm}

\begin{theorem}\label{Thm:FS_Step_Identity}
The intermediary forward substitution vectors $\yvc{k}$ and $\hat{\yv}^{(k)}$ obtained at the conclusion of the $k$th iteration of REF forward substitution on $\vv$ using $L$ and $\hat L$, respectively, are equal,  for $k=1,\dots,n\mi1$.
\end{theorem}
\bproof Upon completion of the $k$th iteration of REF forward substitution, the individual elements of $\yvc{k}$ are connected to IPGE entries according to the equation \citep{esc15rou}:
\begin{equation}
  \yc{k}{i}=\begin{cases}
  \ac{i\mi1}{i}{n+1} & \text{if }  i\le k\\
  \ac{k}{i}{n+1} & \text{if }  i> k;
    \end{cases}
\end{equation}
where column $n\ma1$ of $A$ denotes the right-hand vector on which IPGE is applied ($\vv$ in this theorem). Based on this connection, each element $i>k\ge1$ of $\yvc{k}$ can be equivalently obtained as
\[\yc{k}{i}=\left(L(A^{[k+1]}_{[k],i})D(A^{[k+1]}_{[k],i})^{-1}\right)^{-1}\vv_{[k],i},\]
where $L(A^{[k+1]}_{[k],i})$ and $D(A^{[k+1]}_{[k],i})^{-1}$ are the lower-triangular and diagonal matrices corresponding to the REF-LU factorization of the (nonsingular) submatrix induced by columns $[k+1]=\{1,\dots,k\ma1\}$ and rows $[k]\cup i=\{1,\dots,k,i\}$ of $A$. The above equation can be understood from the observations that element $k\ma1$ of $\yv$ does not change after iteration $k$ of the algorithm and that, if any row $i>k$ of $A$ and $\yv$ is swapped with row $k\ma1$, the new element $k\ma1$ of $\yv$ at iteration $k$ can be obtained using the inverse of the REF lower-triangular and diagonal matrices of the corresponding submatrix of $A$. Furthermore, from the basic identity given by Proposition \ref{prop:nesting},
\begin{equation}
  \left(L(A^{[k+1]}_{[k],i}) D(A^{[k+1]}_{[k],i})^{\mi1}\right)^{-1}=\left((LD^{-1})^{-1}\right)^{[k+1]}_{[k],i},
\end{equation}
where $1\le k< i$. Piecing this together with the above analysis, each entry of $\yvc{k}$ is given by
\begin{equation}
  \yc{k}{i}=\begin{cases}
  \left(\ad{A}^{[i]}_{[i]}\right)_{i}\vv_{[i]} & \text{if }  i\le k\mmV\\
  \left(\ad{A}^{[k+1]}_{[k],i}\right)_{k+1}\vv_{[k],i} & \text{if } i> k.
    \end{cases}
\end{equation}
Since each $\yc{k}{i}$ ($\hat y^{(k)}_i$, resp.) is the product of the bottom row of the adjoint matrix of a submatrix of $A$ ($\hA$, resp.) and a subvector of $\vv$, the proof is completed based on a similar reasoning as Theorem \ref{thm:Inverse_Identity}. \eproof\llV

\begin{theorem}\label{Thm:FS_Step_Identity_b}
The intermediary forward substitution vectors $\zvc{k}$ and $\hat{\zv}^{(k)}$ obtained at the conclusion of the $k$th iteration of REF forward substitution on $\wv$ using $U^T$ and $\hat U^T$, respectively, are equal, for $k=1,\dots,n\mi1$.
\end{theorem}
\bproof The result is obtained through the combined logic of Theorems \ref{thm:Inverse_Identity_b} and \ref{Thm:FS_Step_Identity}. \eproof

\section{REF Rank-One Update Algorithm}\label{Sec:Main_Alg}
This sections introduces and proves the correctness of the featured algorithms and is organized as follows. Section \ref{Subsec:Standard_Version} derives a standard version of the REF rank-one update algorithm, which relies on certain assumptions, and Section \ref{Subsec:Walkthrough} walks through a numerical example. Then, Section \ref{Subsec:Zeros} derives adjustments for dealing with cases when the standard algorithm fails (i.e., when its assumptions do not hold).

\subsection{Standard Version}\label{Subsec:Standard_Version}
The pseudocode of the REF rank-one update algorithm, denoted succinctly as REF-ROU, is given in Algorithm \ref{Rank-one_Alg}. Its steps and their validity are further described in the accompanying proof of correctness. The overall computational complexity of the algorithm is also formally established by describing the number of operations it requires in detail. Beforehand, it is necessary to prove the ensuing lemma.

\begin{lemma}\label{lem:exact_div}
Let $L$ be the lower-triangular matrix factor of the REF-LU factorization of nonsingular matrix $A$, and let $\yvc{k}$ denote the vector output after the $k$th iteration of the stepwise recursion of REF forward substitution performed on a vector $\vv$ using $L$. The following division is exact:
\[\frac{l_{k,k} \yc{k\mi1}{i} - l_{k\mi1,k\mi1} \yc{k}{i}}{\yc{k\mi1}{k}},\]
where $2\le k<i\le n$.
\end{lemma}

\bproof Entry $l_{i,k}\in\Z^1$ can be expanded as
\begin{flalign}
  l_{i,k} & = \frac{l_{i,k}\yc{k\mi1}{k}}{\yc{k\mi1}{k}}\label{R1_Lemma_a}\\
   & = \frac{l_{i,k}\yc{k\mi1}{k}+l_{k,k}\yc{k\mi1}{i}-l_{k,k}\yc{k\mi1}{i}}{\yc{k\mi1}{k}}\label{R1_Lemma_b}\\
   & = \frac{l_{k,k}\yc{k\mi1}{i}-\left(l_{k,k}\yc{k\mi1}{i}-l_{i,k}\yc{k\mi1}{k}\right)}{\yc{k\mi1}{k}}\label{R1_Lemma_c}\\
   & = \frac{l_{k,k}\yc{k\mi1}{i}-l_{k\mi1,k\mi1}\yc{k}{i}}{\yc{k\mi1}{k}}.\label{R1_Lemma_d}
\end{flalign}
The division in the right-hand side of Equation \eqref{R1_Lemma_a} is clearly exact. Its result matches Equations  \eqref{R1_Lemma_b} and  \eqref{R1_Lemma_c}, since in the former the terms added to the numerator cancel out, and in the latter the three numerator terms are only reorganized.  Equation \eqref{R1_Lemma_d} results from the stepwise recursion of REF forward substitution (see Algorithm \ref{FS_Step}); namely, from the formula $\yc{k}{i}=\left(l_{k,k}\yc{k\mi1}{i}-l_{i,k}\yc{k\mi1}{k}\right)/l_{k\mi1,k\mi1}$, the expression within the parenthesis of Equation \eqref{R1_Lemma_c} is substituted with $l_{k\mi1,k\mi1}\yc{k}{j}$, which is a product of two integer entries. Based on the above series of equations, the Equation \eqref{R1_Lemma_d} numerator is a multiple of the denominator. \eproof

\begin{algorithm}[h!]
\SetKwInOut{Input}{input}\SetKwInOut{Output}{output}\SetKw{KwRet}{return}
\Input{$L$, $U$, $\vv$, $\wv$, $Diag(A)$}
{\bf let} $\hat L, \hat U\in\Z^{n\times n}$\\
     $\hat{l}_{1,1}, \hat{u}_{1,1} = a_{1,1} + v_1w_1$\\
    \For{ i = 2,\dots,n}
    {
        $\hat{u}_{1,i} = u_{1,i} + v_1w_i$
        \\
        $\hat{l}_{i,1} = l_{i,1} + v_iw_1$
        \\
        $\hat{l}_{i,i},\hat{u}_{i,i}  = a_{i,i} + v_iw_i$
    }
     $\yvc{1}=$ REF-FS-Step\:($\yvc{0}=\vv, L^1_{[n]}, 1, 1$)            \\
    $\zvc{1}=$ REF-FS-Step\:($\zvc{0}=\wv, (U^T)^1_{[n]}, 1, 1$)
    \\
\For{$k = 2,\dots,n\mi1$}
{
    $\yvc{k} = $ REF-FS-Step\:($\yvc{k\mi1}$, $L^{k}_{[n]}$, $l_{k\mi1,k\mi1}$, $k$)
                \\
    $\zvc{k} = $ REF-FS-Step\:($\zvc{k\mi1}$, $(U^T)^{k}_{[n]}$, $u_{k\mi1,k\mi1}$, $k$)
    \\
    \For{ i = k,\dots,n}
    {
     $\hat{l}_{i,i},\hat{u}_{i,i}  \leftarrow [( \hat{l}_{k\mi1,k\mi1}\hat{l}_{i,i}) - (\hat{u}_{k\mi1,i}\hat{l}_{i,k\mi1})]/\hat{l}_{k\mi2,k\mi2}$

        \If{$i > k$}
        {
        $\hat{l}_{i,k} = [(\hat{l}_{k,k} \yc{k\mi1}{i}) - (\hat{l}_{k\mi1,k\mi1} \yc{k}{i})]/\yc{k\mi1}{k}$

        $\hat{u}_{k,i}= [(\hat{u}_{k,k} \zc{k\mi1}{i}) - (\hat{u}_{k\mi1,k\mi1} \zc{k}{i})]/\zc{k\mi1}{k}$
        }
    }
}
$\hat{l}_{n,n},\hat{u}_{n,n} \leftarrow [( \hat{l}_{n\mi1,n\mi1}\hat{l}_{n,n}) - (\hat{u}_{n\mi1,n}\hat{l}_{n,n\mi1})]/\hat{l}_{n\mi2,n\mi2}$\\
\KwRet{$\hat L,\hat U$}
\caption{REF Rank-one Update Algorithm (REF-ROU)\label{Rank-one_Alg}}
\end{algorithm}

\begin{theorem}\label{thm:Rank-one_Correctness}
Define matrices $A$ and $\hA$ and update vectors $\vv$ and $\wv$ as in Assumption \ref{assu:AMatrix}, and let REF-LU($A)=LD^{-1}U$ and REF-LU($\hA)=\hat L\hat D^{-1}\hat U$. Additionally, let $\yvc{k}$ and $\zvc{k}$ denote the vectors output after the $k$th iteration of the stepwise recursion of REF forward substitution performed on $\vv$ using $L$ and on $\wv$ using $U^T$, respectively, for $k=1,\dots,n\mi1$. Under the assumption that $\yc{k\mi1}{k}\ne0$ and $\zc{k\mi1}{k}\ne0$, for $k=2,\dots,n\mi1$, Algorithm \ref{Rank-one_Alg} (i.e., \textnormal{the standard REF-ROU algorithm}) successfully obtains REF-LU($\hA$) from REF-LU($A$) without roundoff errors.
\end{theorem}
\bproof First, we explain how to obtain the off-diagonal entries of REF-LU($\hA$), whose calculation requires the availability of the updated factorization pivots (i.e., $\hat\rho^k:=\hat l_{k,k}=\hat u_{k,k}$); afterward, we explain how the diagonal of REF-LU($\hA$) is iteratively constructed along the update process to furnish the required pivots.

Since the entries along the first row and column of a REF-LU factorization match the respective input matrix entries, the first row and column of REF-LU($\hA$) are obtained by simply adding $v_1w_i$ to $u_{1,i}=\ac{0}{1}{i}$ and $v_iw_1$ to $l_{i,1}=\ac{0}{i}{1}$, for $i=2,\dots,n$. This is accomplished with the initial for-loop (its first and second statements) without roundoff errors. The bulk of the off-diagonal entries are obtained by utilizing the theoretical insights derived in Section \ref{Sec:Algebra}. In particular, Theorem \ref{Thm:FS_Step_Identity} establishes that $\yvc{k}$ can be equivalently obtained from the $k$th step of REF forward substitution on $\vv$ using either $L$ or $\hat L$ as the lower-triangular matrix factor; similarly, Theorem \ref{Thm:FS_Step_Identity_b} establishes that $\zvc{k}$ can be equivalently obtained from the $k$th step of REF forward substitution on $\wv$ using either $U^T$ or $\hat U^T$. The ensuing arguments leverage the first of these two insights to reconstruct the entries of $\hat L$ from the outputs of REF forward substitution using $L$; parallel arguments are used to reconstruct the entries of $\hat U^T$ from the outputs of REF forward substitution using $U^T$, but they are omitted for brevity. From the stepwise recursion of  REF forward substitution, each entry $\yc{k}{i}$ is obtainable using $\hat L$ as
\[\yc{k}{i}=\frac{\hat l_{k,k}\yc{k\mi1}{i}-\hat l_{i,k}\yc{k\mi1}{k}}{\hat l_{k\mi1,k\mi1}},\]
where $i > k$ and $\hat l_{0,0}=1$. Reorganizing this expression to isolate $\hat l_{i,k}$ gives
\[\hat l_{i,k}=\frac{\hat l_{k,k}\yc{k\mi1}{i}-\hat l_{k\mi1,k\mi1}\yc{k}{i}}{\yc{k\mi1}{k}}.\]
In words, for each $k=2,\dots,n\mi1$, column $k$ of $\hat L$ is reconstructed from the iteration $k\mi1$ and $k$ REF forward substitution vectors---$\yvc{k\mi1}$ and $\yvc{k}$, obtained using $L$---and from updated pivots---$\hat l_{k\mi1,k\mi1} $ and $\hat l_{k,k}$ (see next paragraph); these four required components are calculated (or already available) at the start of the $k$th iteration of the outer for-loop of the algorithm. Because $\yc{k\mi1}{k}\ne0$, from the given assumption, and the above division is exact, according to Lemma \ref{lem:exact_div}, $\hat l_{i,k}$ is successfully obtained free of roundoff error.

REF-ROU obtains the updated factorization pivots as follows. First, Diag($\hA$) is obtained directly from Diag($A$) in the initial for-loop (its third statement). This opening step provides the correct and final value of the first updated pivot ($\hat\rho^1=\hat l_{1,1}=\hat u_{1,1}$), and it populates starting values for diagonal elements $i=2,\dots,n$, which will be finalized one by one in subsequent iterations. In the next for-loop (which starts at $k=2$), the initial step of the inner for-loop provides the correct and final value of the second updated diagonal element ($\hat\rho^2=\hat l_{2,2}=\hat u_{2,2}$); this is established by the fact that the right-hand side of its first statement can be written in terms of IPGE entries as
\[\frac{\hat{l}_{1,1}\hat{l}_{2,2} - \hat{u}_{1,2}\hat{l}_{2,1}}{\hat{l}_{0,0}}=\frac{\hac011\hac022 -  \hac012\hac021}{\hac{-1}00}=\hac122=\hat\rho^2.\]
Similarly, upon completion of the first outer for-loop iteration, the $i$th diagonal entry is equivalent to entry $\hac{1}{i}{i}$, where $\cRan{3}{i}{n}$, meaning these entries have been advanced from iteration 0 to iteration 1 of IPGE and are not yet in their final form. Continuing with this process, upon completion of outer for-loop iteration $k$, diagonal entry $\hat l_{ii}$ with $i>k$ is equivalent to $\hac{k\mi1}{i}{i}$; hence, $\hat\rho^k=\hat l_{k,k}=\hat u_{k,k}$ is finalized during this iteration. The entries required to obtain this IPGE entry are available to the algorithm since
\[\hac{k\mi1}{i}{i}=\frac{\hac{k\mi2}{k\mi1}{k\mi1}\hac{k\mi2}{i}{i} - \hac{k\mi2}{k\mi1}{i}\hac{k\mi2}{i}{k\mi1}}{\hac{k\mi3}{k\mi2}{k\mi2}}=\frac{\hat l_{k\mi1,k\mi1}\hac{k\mi2}ii - \hat u_{k\mi1,i}\hat l_{i,k\mi1}}{\hat l_{k\mi2,k\mi2}},\]
where the second equation is obtained from the relationship between REF-LU$(\hA)$ and the iterative entries of the IPGE algorithm. That is, the entries involved in the calculation of $\hac{k\mi1}{i}{i}$ are drawn from the updated pivots of the previous two iterations and from column $k\mi1$ of $\hat L$ and row $k\mi1$ of $\hat U$, both obtained during iteration $k\mi1$ of the outer for-loop (along with $\hac{k\mi2}{i}{i}$). Since the final updated pivot ($\rho^{(n)}=\hat{l}_{n,n}=\hat{u}_{n,n}$) is calculated using a similar string of operations culminating in the last line of the algorithm, the updated pivots $\hat\rho^k=\hac{k\mi1}{k}{k}$ are calculated correctly using roundoff-error free operations, for $k=1,\dots,n$. \eproof

\begin{theorem}\label{thm:Rank-one_Operations}
REF-ROU requires $O(n^2)$ operations.
\end{theorem}
\bproof Without loss of generality,  we adopt the convention that each addition, subtraction, multiplication, or division of two vector/matrix entries is considered as one operation. The initial for-loop of Algorithm \ref{Rank-one_Alg} requires $6(n\mi1)$ operations, that is, 1 multiplication and 1 addition, for each entry along each of (i) the first row of $\hat U$, (ii) the first column of $\hat L$, and (iii) the diagonal of $\hat L$ (excepting the first element, which is obtained in the previous step). Calculating $\yvc{1}$ and $\zvc{1}$ via Algorithm \ref{FS_Step} entails another $6(n\mi1)$ operations, that is, 2 multiplications and 1 subtraction, for each entry in these vectors (excepting the first element). Executing the $k$th iteration of the second for-loop, where $\cRan{2}{k}{n\mi1}$, requires $20(n\mi k)$ operations, that is, 2 multiplications, 1 addition, and 1 division, for each entry $i>k$ along each of (i) $\yvc{k}$, (ii) $\zvc{k}$, (iii) the $k$th row of $\hat U$, (iv) the $k$th column of $\hat L$, and (v) the diagonal of $\hat L$. In summary, the full first for-loop and each iteration of the second for-loop require $O(n)$ operations; all remaining steps require constant time. Therefore, the full rank-one update requires $O(n^2)$ operations. \eproof

The ensuing corollary provides the worst-case computational complexity of the REF rank-one update  algorithm (abbreviated as WCC(REF-ROU)). The analysis combines the algorithm's $O(n^2)$ required operations and the cost of multiplying/dividing two integers with bit-length $\wb
\le\lceil n\log(\hat{\sigma}\sqrt{n})\rceil$, where $\hat \sigma:=\underset{i,j}{\max\mH}\left\{\max\left(|\ac{0}{i}{j}|,|\hac{0}{i}{j}|\right)\right\}$. The latter costs must be accounted because the REF algorithms entail working with matrix entries with non-fixed precision.

\begin{corollary}\label{cor:full_complexity}
The worst-case computational complexity (WCC) of REF-ROU is given by:
\end{corollary}
\begin{flalign}
\text{WCC(REF-ROU)}&=O(n^2(\wb\log\wb\log\log\wb))\label{eqn:wcc5}\\
 &=O\left(n^3\max(\log^2n\log\log n,\log^2\hat\sigma\log\log\hat\sigma)\right)\label{eqn:wcc6}.
\end{flalign}

It is worth adding that the bound on $\wb$ is somewhat pessimistic \citep{abb01tig,coo11sol}, meaning that the algorithms perform more efficiently in practice than this theoretical measure suggests. Furthermore, it is reasonable to expect efficiency gains when the algorithms are adapted for sparse and other well structured matrices occurring in real-world applications \citep{lou19exa}.

\subsection{Numerical Example}\label{Subsec:Walkthrough}
This subsection provides a numerical application of the REF-ROU algorithm. To start, consider the following input matrix $A\in\Z^{4\times4}$, its REF LU factorization, and update vectors $\vv,\wv\in\Z^{4}$:

\[\small A=\left[\begin{array}{cccc}
    3 & 8 & 7 & 1 \\
    5 & 3 & 5 & 4 \\
    6 & -2 & 1 & 7 \\
    7 & -2 & -6 & 11
\end{array}\right]
,\mmH
\T{REF-LU}(A) =\left[\begin{array}{rrrr}
    3 & 8 & 7 & 1 \\
    5 & -31 & -20 & 7 \\
    6 & -54 & 43 & -29 \\
    7 & -62 & 279 & -89
\end{array}\right]
,\mmH
\vv = \left[\begin{array}{c}
    1 \\
    5 \\
    7 \\
    2
\end{array}\right]
,\mmH
\wv = \left[\begin{array}{c}
    2 \\
    6 \\
    3 \\
    4
\end{array}\right].
\]
REF-LU($A$) above is displayed in a succinct form that merges together the $L$ and $U$ matrices; this is possible because the diagonals of both matrices are identical and the algorithm does not have need for the $D$ matrix. From these inputs, the rank-one matrix $\vv\wv^T$, the update matrix $\hA$, and the REF LU factorization of $\hA$ (i.e., the desired output from REF-ROU) are as follows:
\[\small\vv\wv^T=\left[\begin{array}{rrrr}
    2 & 6 & 3 & 4 \\
    10 & 30 & 15 & 20 \\
    14 & 42 & 21 & 28 \\
    4 & 12 & 6 & 8
\end{array}\right]
,\mmH
\hA=\left[\begin{array}{rrrr}
    5 & 14 & 10 & 5 \\
    15 & 33 & 20 & 24 \\
    20 & 40 & 22 & 35 \\
    11 & 10 & 0 & 19
\end{array}\right]
,\mmH
\T{REF-LU}(\hA) =\left[\begin{array}{rrrr}
    5 & 14 & 10 & 5 \\
    15 & -45 & -50 & 45 \\
    20 & -80 & 10 & 45 \\
    11 & -104 & -50 & -178
\end{array}\right].\]

Next, we describe how to obtain REF-LU($\hA$) by updating REF-LU($A$). The walk-through is divided into four parts: initial steps, outer for-loop iteration $k =2$, outer for-loop iteration $k=3$, and final step. For convenience, the  entries of REF-LU($\hA$) that are finalized after each part are bolded and colored in blue.

\noindent\underline{Initial Steps}. Begin by constructing the entries along the first row and column of REF-LU($\hA$), by taking the matching entries of $A$ and adding the corresponding product of entries from $\vv$ and $\wv$. The specific calculations are as follows:\llV
\[\small\begin{array}{lll}
 \hat{u}_{12} = a_{12} + v_1w_2 \mmmH\mmmH\mmmH \N & \hat{u}_{13} = a_{13} + v_1w_3 \mmmH\mmmH\mmmH \N & \hat{u}_{14} = a_{14} + v_1w_4 \\
   \phantom{\hat{u}_{12}} = 8 + 1(6) = 14 &  \phantom{\hat{u}_{13}} = 7 + 1(3) = 10 &  \phantom{\hat{u}_{14}} = 1 + 1(4) = 5\\\hline
 \hat{l}_{21} = a_{21} + v_2w_1 \mmmH\mmmH\mmmH \N & \hat{l}_{31} = a_{31} + v_3w_1 \mmmH\mmmH\mmmH \N & \hat{l}_{41} = a_{41} + v_4w_1 \\
  \phantom{\hat{l}_{21}} = 5 + 5(2) = 15 &  \phantom{\hat{l}_{31}} = 6 + 7(2) = 20 &  \phantom{\hat{l}_{41}} = 7 + 2(2) = 11.
\end{array}\]
Perform similar operations to obtain the initial elements of the working diagonal of $\hA$:\llV
\[\small\begin{array}{llll}
 \hat{l}_{11} = a_{11} + v_1w_1 \mmmH\mmmH \N & \hat{l}_{22} = a_{22} + v_2w_2 \mmmH\mmmH \N & \hat{l}_{33} = a_{33} + v_3w_3 \mmmH\mmmH \N & \hat{l}_{44} = a_{44} + v_4w_4 \\
  \phantom{\hat{l}_{11}} = 3 + 1(2) = 5 &  \phantom{\hat{l}_{22}} = 3 + 5(6) = 33 &  \phantom{\hat{l}_{33}} = 1 + 7(3) = 22 &  \phantom{\hat{l}_{44}} = 11 + 2(4) = 19.
\end{array}\]
Next, calculate the forward substitution vectors for iteration $k=1$, $\yvc{1}$ and $\zvc{1}$ using Algorithm \ref{FS_Step} and vectors $\yvc{0}=\vv$ and $\zvc{0}=\vv$, respectively. At the end of these steps, these two vectors and the working factorization, denoted as $\hat{L}\backslash\hat{U}$, are given by:

\[\small\yvc{1} = \left[\begin{array}{c}
    1 \\
    10 \\
    15 \\
    -1
\end{array}\right]
,\mmH
\zvc{1} = \left[\begin{array}{c}
    2 \\
    2 \\
    -5 \\
    10
\end{array}\right]
,\mmH
\hat{L}\backslash\hat{U} =\left[\begin{array}{rrrr}
    \blu\bf 5 & \blu\bf14 & \blu\bf10 & \blu\bf5 \\
    \blu\bf15 & 33 & \cdot & \cdot \\
    \blu\bf20 & \cdot & 22 & \cdot \\
    \blu\bf11 & \cdot & \cdot & 19
\end{array}\right].\]

\noindent\underline{Outer for-loop iteration $k=2$}. First, calculate the  forward substitution vectors for iteration $k=2$, $\yvc{2}$ and $\zvc{2}$, using Algorithm \ref{FS_Step} and vectors $\yvc{1}$ and $\zvc{1}$, respectively. The resulting vectors are given by:
\[\small\yvc{2} = \left[\begin{array}{c}
    1 \\
    10 \\
    25 \\
    217
\end{array}\right]
,\mmH
\zvc{2} = \left[\begin{array}{c}
    2 \\
    2 \\
    65 \\
    -108
\end{array}\right].
\]
Second, advance diagonal elements $\hat l_{22}$, $\hat l_{33}$, $\hat l_{44}$ through IPGE pivoting operations:
\[\small\begin{array}{lll}
  \hat{l}_{22} \leftarrow (\hat{l}_{11}\hat{l}_{22} - \hat{u}_{12}\hat{l}_{21})/\hat{l}_{00} \mmmH\mmmH\mmmH \N & \hat{l}_{33} \leftarrow (\hat{l}_{11}\hat{l}_{33} - \hat{u}_{13}\hat{l}_{31})/\hat{l}_{00} \mmmH\mmmH\mmmH \N & \hat{l}_{44} \leftarrow (\hat{l}_{11}\hat{l}_{44} - \hat{u}_{14}\hat{l}_{41})/\hat{l}_{00} \\
  \phantom{\hat{l}_{22}} = [5(33) - 14(15)]/1 = -45 &  \phantom{\hat{l}_{33}} = [5(22) - 10(20)]/1 = -90 &  \phantom{\hat{l}_{44}} = [5(19) - 5(11)]/1 = 40.
\end{array}\]
Third, obtain the off-diagonal entries of the second row of $\hat{U}$ and second column of $\hat{L}$ through the operations:

\[\small\begin{array}{llll}
  \hat{u}_{23} = (\hat{u}_{22}\zc{1}{3} - \hat{u}_{11}\zc{2}{3})/\zc{1}{2} \mmmH\mmmH\mmmH \N & \hat{u}_{24} = (\hat{u}_{22}\zc{1}{4} - \hat{u}_{11}\zc{2}{4})/\zc{1}{2}  \\
  \phantom{\hat{u}_{23}}=[-45(-5) - 5(65)]/2= -50 &  \phantom{\hat{u}_{24}} =[-45(10) - 5(-108)]/2 = 45\\\hline
  \hat{l}_{32} = (\hat{l}_{22}\yc{1}{3} - \hat{l}_{11}\yc{2}{3})/\yc{1}{2} \mmmH\mmmH\mmmH \N & \hat{l}_{42} = (\hat{l}_{22}\yc{1}{4} - \hat{l}_{11}\yc{2}{4})/\yc{1}{2}  \\
  \phantom{\hat{l}_{11}}=[-45(15) - 5(25)]/10 = -80 &  \phantom{\hat{l}_{21}} =[-45(-1) - 5(217)]/10 = -104.
\end{array}\]
At the end of these steps, the working matrix is given by:
\[\small\hat{L}\backslash\hat{U} =\left[\begin{array}{rrrr}
    5 & 14 & 10 & 5 \\
    15 & \blu\bf -45 & \blu\bf -50 & \blu\bf 45 \\
    20 & \blu\bf -80 & -90 & \cdot \\
    11 & \blu\bf -104 & \cdot & 40 \\
\end{array}\right].\]

\noindent\underline{Outer for-loop iteration $k=3$}. First, calculate the forward substitution vectors, $\yvc{3}$ and $\zvc{3}$, using Algorithm \ref{FS_Step} and vectors $\yvc{2}$ and $\zvc{2}$, respectively. These resulting vectors are given by:
\[\small\yvc{3} = \left[\begin{array}{c}
    1 \\
    10 \\
    25 \\
    -76
\end{array}\right]
,\mmH
\zvc{3} = \left[\begin{array}{c}
    2 \\
    2 \\
    65 \\
    89
\end{array}\right].
\]
Second, advance diagonal elements $\hat l_{33}$ and $\hat l_{44}$ through IPGE pivoting operations:
\[\small\begin{array}{lll}
   \hat{l}_{33} \leftarrow (\hat{l}_{22}\hat{l}_{33} - \hat{u}_{23}\hat{l}_{32})/\hat{l}_{11} \mmmH\mmmH\mmmH \N & \hat{l}_{44}
   \leftarrow (\hat{l}_{22}\hat{l}_{44} - \hat{u}_{24}\hat{l}_{42})/\hat{l}_{11} \\
  \phantom{\hat{l}_{33}} = [-45(-90) + 50(-80)]/5 = 10 \mmmH\mmmH\mmmH \N&  \phantom{\hat{l}_{44}} = [-45(40) - 45(-104)]/5 = 576.
\end{array}\]
Third, obtain the off-diagonal entries of the second row of $\hat{U}$ and second column of $\hat{L}$ through the operations:
\[\small\begin{array}{llll}
  \hat{u}_{34} = (\hat{u}_{33}\zc{2}{4} - \hat{u}_{22}\zc{3}{4})/\zc{2}{3} \\
  \phantom{\hat{u}_{34}}=[10(-108) + 45(89)]/65 = 45\\\hline
  \hat{l}_{43} = (\hat{l}_{33}\yc{2}{4} - \hat{l}_{22}\yc{3}{4})/\yc{2}{3} \\
  \phantom{\hat{l}_{43}}=[10(217) + 45(-76)]/25 = -50.
\end{array}\]
At the end of these steps, the working matrix is given by:
\[\small\hat{L}\backslash\hat{U} =\left[\begin{array}{rrrr}
    5 & 14 & 10 & 5 \\
    15 & -45 & -50 & 45 \\
    20 & -80 & \blu\bf10 & \blu\bf45 \\
    11 & -104 & \blu\bf-50 & 576
\end{array}\right].\]
This completes all iterations of the outer for-loop.

\noindent\underline{Final step}. Finalize diagonal element $\hat l_{44}$ through an IPGE pivoting operation:

\[\small\begin{array}{l}
   \hat{l}_{44} = (\hat{l}_{33}\hat{l}_{44} - \hat{u}_{34}\hat{l}_{43})/\hat{l}_{22} \\
   \phantom{\hat{l}_{44}} = [10(576) - 45(-50)]/-45 = \blu\bf-178.
\end{array}\]
REF-LU$(\hA)$ is completed after inserting this updated diagonal element into the previous working matrix.

\subsection{REF-ROU Special Cases and Adjustments}\label{Subsec:Zeros}
In its standard form, REF-ROU fails when the division involved in the calculation of the off-diagonal factorization entries is undefined (see  equations under the conditional statement of Algorithm \ref{Rank-one_Alg}), which is linked to a zero occurring in a specific element of the iterative update vectors. This subsection introduces adjustments for dealing with two related cases.

\underline{\bf Special Case 1}: Zeros occur before any steps of the algorithm are performed. The problem occurs when the update vectors $\vv$ and/or $\wv$ contain a leading sequence of zeros. Note that, when zeros occur in these initial vectors after a leading sequence of nonzeros, the rank-one update algorithm tends to turn these entries into nonzeros as it progresses; however, if it does not, Special Case 2 can be applied.

The adjustment consists of two main parts. To describe them, let $\theta_v$ and $\theta_w$  be the respective indices of the last zero occurring in an uninterrupted sequence extending from the initial elements of $\vv$ and $\wv$, that is,
\[\theta_v=\underset{0\le i\le n}{\max}\{i: v_j=0 \mmH\forall j\le i\} \mmmH \theta_w=\underset{0\le i\le n}{\max}\{i: w_j=0 \mmH\forall j\le i\};\]
where $v_0=w_0=0$ by convention, so that $\theta_v=0$ if $v_1\ne0$ ($\theta_w=0$ if $w_1\ne0$, resp.). In the first part, all entries of REF-LU($A$) along rows $[\theta_v]$  and columns $[\theta_w]$ are copied to REF-LU($\hA$), and the statements required to retrieve those entries in Algorithm \ref{Rank-one_Alg} are bypassed---specifically, the assignment statements (i.e., with ``$=$'') of $\hat{l}_{i,j}$ and $\zvc{k}$, for $j,k\le \theta_w$ and $i>j$; and of $\hat{u}_{i,j}$ and $\yvc{k}$, for $i,k\le \theta_v$ and $j>i$. In the second part of the algorithm, $\yc{\theta_v}{i}$ and $\zc{\theta_w}{j}$ are set to $\yc{0}{i}$ and $\zc{0}{j}$, respectively, for $i>\theta_v$, $j>\theta_w$ (to set up the subsequent iteration of the REF-FS-Step algorithm), and the remaining steps of Algorithm \ref{Rank-one_Alg} are executed. Note that it is still necessary to perform all assignment and reassignment statements (i.e., with ``$\leftarrow$'') of $\hat{l}_{i,i}$ in Algorithm \ref{Rank-one_Alg}, for $i>\max\{\theta_v,\theta_w\}$.

\underline{\bf Special Case 2}: Zeros occur during the execution of the algorithm. The problem occurs when the algorithm encounters $\yc{k\mi1}{k}=0$ or $\zc{k\mi1}{k}=0$ at iteration $k\ge2$ of the outer for-loop. When this happens, it is not possible to evaluate the stated formula for entries $\hat l_{i,k}$ or $\hat u_{k,i}$, respectively, for any $i>k$. Zeros occurring in other elements of the update vectors at iteration $k$ do not pose issues.

The adjustment needed to overcome this special case also consists of two main parts; for simplicity, the discussion focuses on $\yc{k\mi1}{k}$, since the other case is handled similarly (i.e., in the transpose sense). First, it is necessary to permute either columns $k\mi1$ and $k$ of REF-LU($A$) or both rows and columns $k\mi1$ and $k$ of REF-LU($A$). The column permutation requires less effort and can be performed whenever $u_{k\mi1,k}\ne0$; the row and column permutation can always be performed, since $u_{k,k}=\rc{k}\ne0$. Second, the effects of the permutation operation are propagated to the remainder of REF-LU($A$) as well as to the working update vectors and factorization. The permutation may also require performing IPGE pivoting operations on the entry that will become the new $k$th pivot of REF-LU$(\hA)$. The pseudocode of this adjustment subroutine is provided by Algorithm \ref{Special_Case_Adj_y} in  Appendix \ref{app:SC2_alg} (the pseudocode for the case with $\zc{k\mi1}{k}=0$ is provided by Algorithm \ref{Special_Case_Adj_z} therein). The proof of the ensuing theorem provides more details on how this adjustment is performed efficiently via roundoff error-free operations defined for LP-related updates in \cite{esc17rou}.

\begin{theorem}\label{thm:Adjustment_Correctness}
The specified adjustments successfully overcome Special Cases 1 and 2 of REF-ROU, and they require $O\left(\max\{(n\mi\theta_v)^2, (n\mi\theta_w)^2, n(n\mi\max\{\theta_v,\theta_w\})\}\right)$ and $O(n^2)$ operations, respectively.
\end{theorem}
\bproof For Special Case 1, we have that $\hat a_{ij}=a_{ij}+v_iw_j=a_{ij}$ whenever $v_i=0$ or $w_i=0$, meaning that all entries along the first $\theta_v$ rows and the first $\theta_w$ columns of $A$ and $\hat A$ are identical. Accordingly, the entries of REF-LU($\hat A)$ along rows $[\theta_v]$ and columns $[\theta_w]$ match those of REF-LU($A$) and can be copied directly, allowing the associated operations to be skipped. Calculating the off-diagonal entries along all remaining rows of $\hat U$ and columns of $\hat L$ requires $O\left((n\mi\theta_v)(n\mi\theta_w)-(n\mi\max\{\theta_v,\theta_w\})\right)$ operations (the expression corresponds to the number of entries in the remaining rectangular submatrix, minus the elements that fall along the diagonal). The algorithm also skips the calculation of $\yvc{k}$, for $k\le \theta_v$, and $\zvc{k}$, for $k\le \theta_w$, owing to the leading zeros in $\vv$ and $\wv$ and applicable shortcuts that can be exploited in the IPGE algorithm \citep{lee95fra}; all remaining iterations of the REF forward substitution algorithm on $\vv$ and $\wv$ require $O(n\mi\theta_v)^2$ and  $O(n\mi\theta_w)^2$ total operations, respectively. Finally, obtaining all diagonal elements $\hat l_{i,i}$, with $i>\max\{\theta_v,\theta_w\}$, requires $O\left(n(n\mi\max\{\theta_v,\theta_w\})\right)$ operations, since each such entry must undergo exactly $i$ steps in the IPGE algorithm. The number of operations required by the reduced forward substitution algorithms may be higher or lower than that required to obtain the diagonal entries, depending on the values of $\theta_v,\theta_w$, and $n$; however, they dominate the number of operations required to obtain the off-diagonal entries. Therefore, overcoming Special Case 1 requires $O\left(\max\{(n\mi\theta_v)^2, (n\mi\theta_w)^2, n(n\mi\max\{\theta_v,\theta_w\})\}\right)$ operations.

For Special Case 2, recall that subvectors $A^{1}_{[k]},A^{2}_{[k]}\dots, A^{k}_{[k]}$ are linearly independent, for all $k\ge1$, based on Assumption \ref{assu:pivot} (if row and/or column permutations are required during the factorization of $A$, these arguments would apply to the permuted matrix, say $P_rAP_c$). This also implies that subvectors $A^{1}_{[i]},A^{2}_{[i]}\dots, A^{k}_{[i]}$ are linearly independent, where $\cRan{k}{i}{n}$. Next, it is useful to re-express $\yc{k\mi1}{k}$ as:\llV
\begin{flalign*}
  \yc{k\mi1}{k} & = \ac{k\mi1}{k}{n+1} \\
   & =\det\left({A}^{[k\mi1],n+1}_{[k]}\right),
\end{flalign*}
where the first equation is obtained from the properties of REF forward substitution and the second from the relationship between IPGE entries and the determinants of submatrices of $A$ (see Equation \eqref{eqn:subdet}). Based on this relationship, $\yc{k\mi1}{k}=0$ implies that $\vv_{[k]}$ is in the span of the linearly independent subvectors $A^1_{[k]}, \dots, A^{k\mi1}_{[k]}$, that is,
\begin{equation}\label{eqn_span_a}
\sum_{j=1}^{k\mi1}\alpha_jA^j_{[k]}=\vv_{[k]},
\end{equation}
for some $\alpha_1,\dots,\alpha_{k\mi1}\in\R^1$. When this occurs, however, $\vv_{[k]}$ cannot be simultaneously in the span of subvectors $A^1_{[k]},\dots, A^{k\mi2}_{[k]},A^{k}_{[k]}$. We show this through contradiction by assuming that $\vv_{[k]}$ can be expressed as:
\begin{equation}\label{eqn_span_b}
\sum_{j=1}^{k-2}\beta_jA^j_{[k]}+\beta_{k\mi1}A^{k}_{[k]}=\vv_{[k]},
\end{equation}
for some $\beta_1,\dots,\beta_{k\mi1}\in\R^1$. Combining \eqref{eqn_span_a} and \eqref{eqn_span_b} gives that:
\begin{flalign*}
 \sum_{j=1}^{k\mi2}\beta_jA^j_{[k]}+\beta_{k\mi1}A^{k}_{[k]} &= \sum_{j=1}^{k\mi1}\alpha_jA^j_{[k]}  \\
  \beta_{k\mi1}A^{k}_{[k]} &=  \sum_{j=1}^{k\mi2}(\alpha_j-\beta_j)A^j_{[k]}+\alpha_{k\mi1}A^{k\mi1}_{[k]}  \\
 A^{k}_{[k]} &=  \sum_{j=1}^{k\mi2}\left(\frac{\alpha_j-\beta_j}{\beta_{k\mi1}}\right)A^j_{[k]}+\left(\frac{\alpha_{k\mi1}}{\beta_{k\mi1}}\right)A^{k\mi1}_{[k]}.
\end{flalign*}
The bottom equation indicates that subvector $A^{k}_{[k]}$ is in the span of subvectors $A^1_{[k]}, \dots, A^{k\mi1}_{[k]}$, contradicting the assumption that subvectors $A^1_{[k]},\dots, A^{k}_{[k]}$ are linearly independent. Therefore, $\vv_{[k]}$ cannot also be in the span of subvectors $A^1_{[k]},\dots, A^{k\mi2}_{[k]},A^{k}_{[k]}$. This implies that, if columns $k\mi1$ and $k$ of REF-LU($A$) are exchanged, the new value of $\yc{k\mi1}{k}$---equal to $\det\left({A}^{[k\mi2],k,n+1}_{[k]}\right)$---must be nonzero. Through a similar line of reasoning, if both rows $k\mi1$ and $k$ and columns $k\mi1$ and $k$ of the original factorization are exchanged, the new value of $\yc{k\mi1}{k}$ must be nonzero since
\[\det\left({A}^{[k\mi2],k,n+1}_{[k-2],k,k\mi1}\right)=-\det\left({A}^{[k\mi2],k,n+1}_{[k]}\right)\ne0.\]
In summary, when $\yc{k\mi1}{k}=0$, permuting columns $k\mi1$ and $k$ or both rows and columns $k\mi1$ and $k$ of REF-LU($A$) changes the value of $\yc{k\mi1}{k}$ to be nonzero.

Finally, the adjustment for Special Case 2 can be performed efficiently and without roundoff error via an Adjacent Pivot Column Permutation (APCP) or an Adjacent Pivot Diagonal Permutation (APDP) of REF-LU($A$), two operations originally defined in \cite{esc17rou} for performing LP-related updates. Both subroutines are  roundoff error-free and require $O(n)$ operations. In the worst case, one of these two subroutines must be called in the initial steps of REF-ROU and prior to the start of each of its $n\mi2$ outer for-loop iterations---this occurs precisely when $\vv$ is a multiple of $A^1_{[n]}$. However, because $\yc{k\mi1}{k}$ is guaranteed to become non-zero after each call, this special case still requires $O(n^2)$ total operations. \eproof

The lower number of operations to overcome Special Case 1 stems from a specific type of sparsity, namely a long sequence of leading zeros in the update vectors. It may be possible to exploit a more general type of sparsity through preliminary permutations on $\vv$ and/or $\wv$ and on the respective rows and/or columns of REF-LU$(A)$. These subroutines widen the applicability of the Special Case 1 adjustments, but they entail additional operations. A detailed analysis of the computational trade-offs is provided in Appendix \ref{app:sparsity}.

\section{Computational Tests}\label{Sec:Tests}
This section presents three experiments to evaluate different computational aspects of the featured algorithms. All instances consist of fully dense, randomly generated matrices whose entries are drawn uniformly from the non-zero integers in the interval [-100,100]. The experiments were carried out on a computer with an Intel(R) Xeon(R) CPU E5-2680 @ 2.40 GHz with 64 GB RAM. The code was written in C++ using the unlimited-precision GNU Multiple Precision Arithmetic Library (GMP).

The first two experiments compare the times required to compute the REF LU factorization of $\hA=A+\vv\wv^T\in\mathbb{Z}^{n\times n}$ from scratch and to obtain the factorization from an existing REF LU factorization of $A\in\mathbb{Z}^{n\times n}$; the tested matrix dimensions are $n=2^4, 2^5,\dots,2^{10}$. The experiments differ based on how the update vectors  $\vv,\wv\in\Z^n$ are initialized. The first generates them randomly using the same specifications as the entries of $A$. The second experiment copies a leading segment of $\vv$ from a column of $A$, and it generates its remaining entries as in the first experiment. Specifically, a column index $c$ is drawn uniformly from the integers in the interval $[1,n]$, and a row parameter $r$ is drawn uniformly from the integers in the interval $[c,n]$; then the leading segment of $\vv$ is set as $\vv_{[r]}=A^c_{[r]}$, and the entries of its trailing segment $\vv_{[n]\backslash[r]}$ are generated randomly. Based on the analysis from  Section \ref{Subsec:Zeros}, this second way of generating $\vv$ is guaranteed to trigger $(r-c)$ calls to the Special Case 2 adjustments.

The results are shown in Table \ref{table:Exps} under columns ``REF-LU($\hA$)'' and ``REF-ROU'', respectively. Therein, for each tested value of $n$, the arithmetic mean and standard deviation of run-times over 30 different instances are reported in seconds (s), rounded to three decimals. For the second experiment, the table also reports summary statistics on the number of calls to the Special Case 2 adjustments (under column ``SC2-Calls''), rounded to one decimal.

\begin{table}[!ht]\small
\vspace{.4in}\centering\def\arraystretch{1}\arraycolsep=3pt \caption{Summary of Computational Results --- Exact Factorization vs. Exact Rank-one Update}\label{table:Exps}
\subfloat[Experiment 1]{
$\begin{array}{|r|rr|rr|}\hline
&\multicolumn{4}{c|}{\small\T{\bf Run-times (s)}}\\
&\multicolumn{2}{c|}{\small\T{\bf REF-LU}(\hA)} &\multicolumn{2}{c|}{\small\T{\bf REF-ROU}}\\
\mH n\N &\small\T{AVG}  &\small\T{ \N SD\N}  &\small\T{AVG}
&\small\T{ \N SD\N} \\\hline
16	&	0.000	&	0.000	&	0.000	&	0.000	\\
32	&	0.001	&	0.000	&	0.000	&	0.000	\\
64	&	0.015	&	0.003	&	0.004	&	0.001	\\
128	&	0.161	&	0.015	&	0.019	&	0.002	\\
256	&	2.484	&	0.016	&	0.150	&	0.001	\\
512	&	57.061	&	1.652	&	1.554	&	0.048	\\
1024	&	1357.219	&	17.664	&	17.763	&	0.114\\
\hline
\end{array}$\label{table:Exp1}}\hspace{.5cm}
\subfloat[Experiment 2]{
$\begin{array}{|r|rr|rr|rr|}\hline
&\multicolumn{4}{c|}{\small\T{\bf Run-times (s)}} &\multicolumn{2}{c|}{}\\
&\multicolumn{2}{c|}{\small\T{\bf REF-LU}(\hA)} &\multicolumn{2}{c|}{\small\T{\bf REF-ROU}}
&\multicolumn{2}{c|}{\small\T{\bf SC-2 Calls}}\\
\mH n\N &\small\T{AVG}  &\small\T{ \N SD\N}  &\small\T{AVG}
&\small\T{ \N SD\N}  &\small\T{AVG} &\small\T{ \N SD\N}\\\hline
16	&	0.000	&	0.000	&	0.000	&	0.000	&	4.6	&	3.8	\\
32	&	0.002	&	0.001	&	0.001	&	0.000	&	7.8	&	6.5	\\
64	&	0.016	&	0.003	&	0.004	&	0.001	&	13.2	&	11.5	\\
128	&	0.161	&	0.015	&	0.023	&	0.006	&	31.8	&	28.7	\\
256	&	2.492	&	0.043	&	0.182	&	0.031	&	52.6	&	49.3	\\
512	&	55.578	&	0.226	&	1.893	&	0.322	&	118.9	&	104.4	\\
1024	&	1354.077	&	9.807	&	22.629	&	3.720	&	281.2	&	221.0\\\hline
\end{array}$\label{table:Exp2}}
\end{table}

\begin{figure}[H]
    \centering
    \hspace*{-10mm}
    \includegraphics[width=4in]{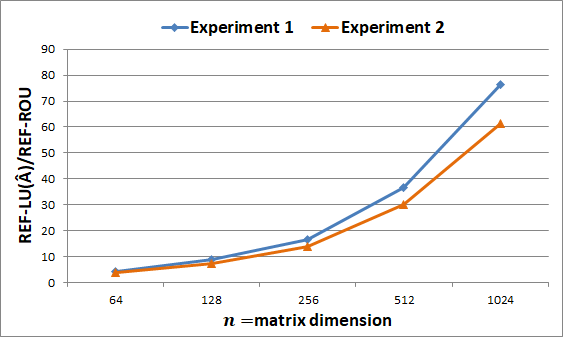}
    \caption{Run-time Ratios of Exact Factorization to the Rank-one Update Algorithm}
    \label{fig:Graph_1}
\end{figure}

As expected, the results demonstrate that obtaining REF-LU($\hA$) by performing a REF rank-one update is orders of magnitude faster than constructing this REF LU factorization from scratch. In the first experiment, the exact update required less than 18 seconds while the refactorization required over 22 \textit{minutes}, on average, for the largest matrix tested---which has over 1,000,000 non-zeros. Stated otherwise, the update could be performed close to 76-times in the time it takes to build the exact factorization for a matrix of this size. In the second experiment, this performance advantage decreased to 61-times, on account of the high number of calls to the ACPU algorithm that are needed to adjust for Special Case 2, based on the way that $\vv$ is generated (i.e., we force $\yc{k\mi1}{k}=0$ for an unrealistically high number of iterations). Even with this extra effort, however, the average REF-ROU times are under 23 seconds. For completeness, Figure \ref{fig:Graph_1} plots the average of the performance ratios REF-LU($\hA$)/REF-ROU, for both experiments and all tested values of $n$.

It is important to add that no calls to ACPU were required for any of the 210 instances tested in Experiment 1, indicating that it is highly unlikely for a leading segment of $\vv$ ($\wv$, resp.) to be linearly dependent on the corresponding segments of the first columns (rows, resp.) of $A$, when the input matrix and update vectors are fully dense and randomly generated as in the featured experiments.

\vspace{-1mm}

The third experiment compares the run-times required to perform a roundoff-error-free column replacement of REF-LU via REF-ROU versus with the push-and-swap approach (P\&S), which was specifically developed for performing simplex algorithm-related updates in \cite{esc17rou} (see Section \ref{Subsec:REF-Up}). To explain this additional use of the rank-one updates, assume that an exiting column $\av^k:=A^{k}_{[n]}$ (the $k$th column of $A$), where $\cRan{1}{k}{n}$, is to be replaced by an incoming column, denoted as $\av^{n+1}\in\mathbb{Z}^n$. This replacement of columns of $A$ can be expressed succinctly as
\begin{equation}
\hA= A + (\av^{n+1}-\av^k)\mH\xlH\ev^T_k,
\label{eqn:ColRep}
\end{equation}
where $\ev_k$ is the $k$th elementary vector of length $n$. That is, Equation \eqref{eqn:ColRep} constitutes a special case of Equation \eqref{eqn:rank-one_SLE} in which $\gamma=1$, $\vv=(\av^{n+1}-\av^k)$, and $\wv=\ev_k$.

The experiment implements the most computationally demanding single-column replacement by setting $\wv=\ev_1$ (i.e., the first column of $A$ is always the replaced column). Another difference from the first two experiments is that the tested matrix dimensions are $n=200, 300,\dots,1000$; this removes the smallest matrices from the previously tested range, for which both approaches finish in tiny fractions of a second, and it helps to better contrast their performance. The results are summarized in Table \ref{table:Exp3}. Therein, the columns under ``REF-ROU'' and ``P\&S'' report, for each approach, the arithmetic mean and standard deviation of run-times obtained over 30 different instances (each with its own underlying REF LU factorization); the columns under ``REF-ROU/P\&S'' report statistics on the performance of the rank-one update approach relative to the push-and-swap approach.

\begin{table}[!ht]\small
\vspace{.4in}\centering\def\arraystretch{1}\arraycolsep=3pt \caption{Summary of Computational Results --- Column Replacement via REF-ROU vs. via P\&S}{
$\begin{array}{|r|rr|rr|rr|}\hline
&\multicolumn{4}{c|}{\small\T{\bf Run-times (s)}} &\multicolumn{2}{c|}{\small\T{\bf REF-ROU/}}\\
&\multicolumn{2}{c|}{\small\T{\bf REF-ROU}} &\multicolumn{2}{c|}{\small\T{\bf P\&S}}
&\multicolumn{2}{c|}{\small\T{\bf P\&S}}\\
\mH n\N &\small\T{AVG}  &\small\T{ \N SD\N}  &\small\T{AVG}
&\small\T{ \N SD\N}  &\small\T{AVG} &\small\T{ \N SD\N}\\\hline
100	&	0.012	&	0.001	&	0.004	&	0.000	&	2.757	&	0.218	\\
200	&	0.079	&	0.001	&	0.037	&	0.001	&	2.121	&	0.034	\\
300	&	0.274	&	0.002	&	0.137	&	0.001	&	2.006	&	0.009	\\
400	&	0.734	&	0.011	&	0.377	&	0.008	&	1.948	&	0.025	\\
500	&	1.594	&	0.025	&	0.822	&	0.005	&	1.940	&	0.029	\\
600	&	3.004	&	0.018	&	1.595	&	0.007	&	1.884	&	0.011	\\
700	&	5.072	&	0.029	&	2.702	&	0.015	&	1.877	&	0.009	\\
800	&	8.269	&	0.057	&	4.510	&	0.050	&	1.834	&	0.015	\\
900	&	12.573	&	0.042	&	6.906	&	0.023	&	1.821	&	0.009	\\
1000	&	18.600	&	0.136	&	10.295	&	0.087	&	1.807	&	0.016	\\
\hline
\end{array}$\label{table:Exp3}}
\end{table}

The computational results demonstrate that REF-ROU is a less efficient approach for performing column replacement updates; it takes approximately 1.8x of the computational times achieved by P\&S on the largest matrices tested, and its performance worsens on smaller matrices. The superior performance of the P\&S approach can be explained by the lower number of operations required in its $k$th iteration, where $\cRan{2}{k}{n\mi1}$ (the first iteration is ignored for simplicity). The $k$th iteration of P\&S entails $12(n\mi k\ma1)$ operations; that is, 2 multiplications, 1 addition or subtraction, and 1 division, for updating each entry $i>k$ along each of (i) the $k$th row of $\hat U$ (which has one more column than $U$ until the update concludes), (ii) the $k$th column of $\hat L$, and (iii) $\yvc{k}$  (forward substitution is performed in full on the incoming column before modifications are made to the factorization; for the purpose of this analysis, it is combined with the other operations). The $k$th iteration of REF-ROU requires $20(n\mi k)$ operations (see the proof to Theorem \ref{thm:Rank-one_Operations}), which is approximately 67\% more than P\&S requires---note that the ratio in computational times seems draw closer to this theoretical value as $n$ increases in Table \ref{table:Exp3}. Besides the higher number of operations that REF-ROU entails, the performance differences can be attributed to the fact that REF-ROU requires a separate $n\times n$ matrix to calculate and store the updated factorization; conversely, P\&S performs the column replacement update on top of the existing REF-LU factorization (augmented by one column appended to the right of $U$).

Altogether, these findings highlight the respective practical advantages of the two update approaches. In particular, although REF-ROU can be also applied to perform LP-related updates of REF-LU in $O(n^2)$ operations, it is outperformed by P\&S in this special case. On the other hand, while P\&S requires fewer operations (although still $O(n^2)$) and approximately half of the memory of REF-ROU to perform a column replacement update, it is inefficient to deploy it to perform a general rank-one update---since this would entail $O(n^3)$ operations, as is discussed in Section \ref{Subsec:REF-Up}.

\section{Conclusion and Future Work}\label{Sec:Conclusion}
This work introduces a direct solution approach for efficiently solving systems of linear equations (SLEs) obtained from rank-one modifications, which are core subroutines used in  nonlinear programming (NLP) and many scientific applications. More specifically, it introduces algorithms for updating existing exact LU and Cholesky factorizations using integer-preserving arithmetic, rather than building new exact factorizations each time the current SLE is modified. The formal guarantees of the algorithms are established through the derivation of theoretical insights, and their computational advantages are supported with computational experiments, which demonstrate upwards of 75x-improvements over exact factorization run-times on fully dense matrices with over one million entries.  Altogether, the exact rank-one updates serve as a foundation for enabling the implementation of the roundoff-error-free (REF) optimization framework, originally developed for linear programming (LP), within broader classes of optimization problems.

The proposed algorithms are likely to be computationally prohibitive on matrices significantly larger than those tested in this work (at least without a corresponding increase in computational resources). To enhance their viability for real-world applications, our future work will seek to develop fully sparse implementations consisting of sparse algorithms and data structures. Analogous efforts in the context of exact LP in \cite{lou19exa} enabled the solution of benchmark instances from the BasisLIB\_INT repository (see \cite{coo11sol}) consisting of SLEs with up to 50,000 rows/columns in relatively reasonable times---in fact, only 3 of 276 these instances  surpassed 1,350 seconds (i.e., the average run-time of the largest instances tested herein). Therefore, it is reasonable to conjecture that the sparse versions of the REF rank-one update algorithms will be orders of magnitude faster than their dense counterparts. As additional avenues to make the REF algorithms more practical for NLP, we will seek to exploit specially structured matrices (e.g., band matrices) occurring in specific applications and to analyze their implications (e.g., tighter bounds on $\beta_{\max}$). A parallel direction is to tailor the proposed algorithms to sparse Cholesky factorizations \citep{dav01mul,dav05row} and related applications (e.g., \cite{her18fac,her20spa}).

It is important to emphasize that the exact methods developed in this work are not intended to replace their inexact counterparts but rather to complement them for the purpose of guaranteeing fast and valid solutions to NLP problems. Along these lines, \cite{web19sol} recently introduced an iterative refinement technique that quickly calculates high-precision KKT solutions of convex and non-convex quadratic programming problems. However, because the guarantees of this approach---and other indirect solution approaches (e.g.,  \citep{gil15met})---depend on certain technical assumptions regarding the conditioning of the inputs, they are not fail-proof. Hence, as another future research direction, we will explore how the REF rank-one update algorithms could be used to supplement this algorithm, analogous to how exact LU factorization is used alongside LP iterative refinement within the \texttt{SoPlex} solver \citep{soplex,wun96par}, which is part of the SCIP Optimization Suite \citep{scip}. More generally, a promising extension for this research is to explore when scaling and other traditional numerical conditioning techniques are sufficient to guarantee the validity of solver outputs and to reserve the REF factorization algorithms for the most numerically challenging instances.

Finally, it is worth mentioning that roundoff errors and their implications may not be an important concern in many situations, especially when the errors inherent in the problem inputs may be larger in magnitude than one would expect from the use of floating-point computations. In the latter cases, obtaining a highly accurate or exact solution may not be justified by the added computational effort required. While this research aims to make the use of exact methods more viable, practitioners and subject-matter experts should weigh these factors to determine the suitable level of precision needed for the problem at hand.

\section{Acknowledgments}\label{Sec:Acknowledgements}
Thanks are extended to Venkata Saisrikar Gudivada for providing one of the initial insights that motivated to this work. This research was supported in part through high performance computing
resources provided by Arizona State University.

\newpage
\begin{APPENDIX}{}
\section{Special Case 2 Adjustment Algorithms}\label{app:SC2_alg}
This appendix introduces two algorithms associated with Special Case 2 of REF-ROU (see Section \ref{Subsec:Zeros}). It also describes a further subcase that may occur during the execution of these adjustments and outlines the respective steps needed to overcome it. Algorithm \ref{Special_Case_Adj_y} provides the necessary steps for ensuring that $\yc{k\mi1}{k}\ne0$, and Algorithm \ref{Special_Case_Adj_z} provides the necessary steps for ensuring that $\zc{k\mi1}{k}\ne0$, where $2\le k\le n\mi1$. The two subroutines are assumed to be embedded within Algorithm \ref{Rank-one_Alg}, specifically, they are executed prior to each call of the REF forward substitution algorithm. The first line of these algorithms checks whether the element of the corresponding iterative REF forward substitution vector that is to be utilized as the divisor in the calculations of the updated factorization's off-diagonal entries will be zero-valued in the next iteration, prior to actually having to perform the current iteration. Expressly, it checks whether
$\yc{k}{k+1}=\hat{l}_{k\mi1,k\mi1}\yc{k\mi1}{k+1}-\hat{l}_{k+1,k\mi1}\yc{k\mi1}{k\mi1}=0$ or $\zc{k}{k+1}=\hat{u}_{k\mi1,k\mi1}\zc{k\mi1}{k+1}-\hat{u}_{k\mi1,k+1}\zc{k\mi1}{k\mi1}=0$. Note that these zero-valued elements would not pose an issue in the current REF-ROU iteration, where they appear in the numerator in the calculations of entries $\hat l_{k+1,k}$ and $\hat u_{k,k+1}$; however, they become problematic in the next iteration, where they become the divisors $\yc{k\mi1}{k}$ and $\zc{k\mi1}{k}$, respectively. This impending issue is avoided in the current iteration by performing the requisite row and/or column permutations of REF-LU($A$) via an Adjacent Pivot Column Permutation update (APCPU), an Adjacent Pivot Row Permutation update (APRPU), or an Adjacent Pivot Diagonal Permutation update (APDPU). The three roundoff-error-free subroutines are defined in \cite{esc17rou}. In addition, the corresponding permutations are performed on the update vectors and working updated factorization and, in the case of APCPU or APRPU, $O(n)$ IPGE pivoting operations  are performed on the entry that will become the new $k$th pivot element of REF-LU$(\hA)$.

\begin{algorithm}[h!]\small
\SetKwInOut{Input}{input}\SetKwInOut{Output}{output}\SetKw{KwRet}{return}
\mmV
{\bf if} $\hat{l}_{k\mi1,k\mi1}\yc{k\mi1}{k+1}-\hat{l}_{k+1,k\mi1}\yc{k\mi1}{k\mi1}=0$\\

\mmmH {\bf if} $\hat u_{k\mi1,k}\ne 0$\\
\mmmH\mmmH Permute columns $k\mi1,k$ of $A$ and $\hat U$\\
\mmmH\mmmH Permute elements $k\mi1,k$ of $\wv$ and $\zvc{k\mi1}$\\
\mmmH\mmmH REF-LU($A$) $\leftarrow$ APCPU(REF-LU($A$), $k\mi1\leftrightarrow k$)\\
\mmmH\mmmH $\backslash*$ Perform $k\mi2$ IPGE pivoting operations on the entry that will become the new $k$th pivot:\\
\mmmH\mmmH $\hat l_{k,k}\leftarrow\ac{0}{k}{k}+v_kw_k$\\
 \mmmH\mmmH \For{ i = 1,\dots,k\mi2}
    {
        \mmmH\mmmH\mmmH $\hat l_{k,k}\leftarrow\hat l_{i,i}\hat l_{k,k}-\hat u_{i,k}\hat l_{k,i}$\\
        \mmmH\mmmH\mmmH {\bf if} $i\ge2$\\
        \mmmH\mmmH\mmmH\mmmH $\hat l_{k,k}\leftarrow \hat l_{k,k}/\hat l_{i\mi1,i\mi1}$\\
  }
  \mmmH\mmmH $\hat u_{k,k}\leftarrow \hat l_{k,k}$\\
  \mmmH {\bf else}\\
    \mmmH\mmmH Permute columns $k\mi1,k$ of $A$ and $\hat U$\\
    \mmmH\mmmH Permute rows $k\mi1,k$ of $A$ and $\hat L$\\
    \mmmH\mmmH Permute elements $k\mi1,k$ of $\vv,\wv,\yvc{k\mi1}$, and $\zvc{k\mi1}$\\
    \mmmH\mmmH REF-LU($A$) $\leftarrow$ APDPU(REF-LU($A$), $k\mi1\leftrightarrow k$)\\

\caption{ROU Special Case 2 Subroutine for Avoiding $\yc{k-1}{k}=0$\label{Special_Case_Adj_y}}
\end{algorithm}

\begin{algorithm}[h!]\small
\SetKwInOut{Input}{input}\SetKwInOut{Output}{output}\SetKw{KwRet}{return}

\mmV
{\bf if}
$\hat{u}_{k\mi1,k\mi1}\zc{k\mi1}{k+1}-\hat{u}_{k\mi1,k+1}\zc{k\mi1}{k\mi1}=0$\\
\mmmH {\bf if} $\hat l_{k,k\mi1}\ne0$\\
\mmmH\mmmH Permute rows $k\mi1,k$ of $A$ and $\hat L$\\
\mmmH\mmmH Permute elements $k\mi1,k$ of $\vv$ and $\yvc{k\mi1}$\\
\mmmH\mmmH REF-LU($A$) $\leftarrow$ APRPU(REF-LU($A$), $k\mi1\leftrightarrow k$)\\
\mmmH\mmmH $\backslash*$ Perform $k\mi2$ IPGE pivoting operations on the entry that will become the new $k$th pivot:\\
\mmmH\mmmH $\hat l_{k,k}\leftarrow\ac{0}{k}{k}+v_kw_k$\\
 \mmmH\mmmH \For{ i = 1,\dots,k\mi2}
    {
        \mmmH\mmmH\mmmH $\hat l_{k,k}\leftarrow\hat l_{i,i}\hat l_{k,k}-\hat u_{i,k}\hat l_{k,i}$\\
        \mmmH\mmmH\mmmH {\bf if} $i\ge2$\\
        \mmmH\mmmH\mmmH\mmmH $\hat l_{k,k}\leftarrow \hat l_{k,k}/\hat l_{i\mi1,i\mi1}$\\
  }
  \mmmH\mmmH $\hat u_{k,k}\leftarrow \hat l_{k,k}$\\
  \mmmH {\bf else}\\
    \mmmH\mmmH Permute columns $k\mi1,k$ of $A$ and $\hat U$\\
    \mmmH\mmmH Permute rows $k\mi1,k$ of $A$ and $\hat L$\\
    \mmmH\mmmH Permute elements $k\mi1,k$ of $\vv,\wv,\yvc{k\mi1}$, and $\zvc{k\mi1}$\\
    \mmmH\mmmH REF-LU($A$) $\leftarrow$ APDPU(REF-LU($A$), $k\mi1\leftrightarrow k$)\\

\caption{ROU Special Case 2 Subroutine for Avoiding $\zc{k-1}{k}=0$\label{Special_Case_Adj_z}}
\end{algorithm}

We discuss the subcase that may occur during the execution of the Special Case 2 adjustments with respect to Algorithm \ref{Special_Case_Adj_y} (modifications in the transpose sense are required for Algorithm \ref{Special_Case_Adj_z}). When the algorithm permutes both rows and column $k\mi1$ and $k$ of REF-LU($A$) in the current iteration, this may yield $\zc{k\mi1}{k}=0$ in the next iteration. That is, $\wv_{[k]}$ would be in the span of the linearly independent subvectors $A^{[k]}_1, \dots, A^{[k]}_{k\mi1}$ following these permutations, thereby complicating the retrieval of the corresponding row of $\hat U$. Since the occurrence of Special Case 2 on dense matrices is already atypical (see Section \ref{Sec:Tests}), only concise descriptions of this even rarer case and its requisite adjustments are provided. To identify its potential occurrence, check that $\hat u_{k,k}\zc{k\mi1}{k+1}-\hat u_{k\mi1,k+1}\zc{k\mi1}{k}=0$ immediately after the else-statement in Algorithm \ref{Special_Case_Adj_y} is called---this mathematical expression corresponds to the REF-FS-Step statement that would be used to obtain the element that will become $\zc{k\mi1}{k}=0$ in the next iteration. If the given conditional statement is true, perform the following steps before continuing with Algorithm \ref{Special_Case_Adj_y}. First, proceed to obtain $\zvc{k}$ and the $k$th row of $\hat U$ using the associated Algorithm 2 steps. Second \textit{backtrack} the entries along row $k$ of $\hat U$, that is, calculate what their value would be in the preceding IPGE iteration ($k\mi2$); the additional entries needed to carry out this subroutine are located along column $k\mi1$ of $\hat L$ and row $k\mi1$ of $\hat U$ \citep{esc17rou}. Third, permute rows $k\mi1,k$ of $\hat U$ and afterwards advance the new row-$k$ entries to the next IPGE iteration ($k\mi1$) using column $k\mi1$ of $\hat L$ as the pivot column and the new row $k\mi1$ of $\hat U$ as the pivot row (see \eqref{eqn:ipge}). The entries of $\zvc{k}$ and rows $k\mi1,k$ of $\hat U$ are finalized at this juncture, meaning that any reassignments of these entries are skipped in the remaining steps of Algorithm \ref{Special_Case_Adj_y} and within the current for-loop iteration of Algorithm 2. Because each of the three above steps requires $O(n)$ operations, these extra adjustments keep REF-ROU to $O(n^2)$ operations.

\section{Exploiting the General Sparsity of the Update Vectors}\label{app:sparsity}
\noindent This appendix examines how to exploit a more general type of sparsity in $\vv$ and/or $\wv$ within REF-ROU; for simplicity, the discussion focuses primarily on the second update vector. When $\theta_w$ is small but $\wv$ is nonetheless relatively sparse, it is possible to apply the Special Case 1 adjustments after performing preliminary permutations of the update vector so as to render it with a longer sequence of leading zeros. For REF-ROU to be properly defined, however, matching permutations of the columns of the existing REF LU factorization are also required. Furthermore, to preserve the integer-preserving properties of the factorization, the permutation of any two columns must be performed via a sequence of calls to APCP (defined in the proof of Theorem \ref{thm:Adjustment_Correctness}). That is, the permutation of columns $k$ and $k'$ of REF-LU$(A)$, where $1\le k< k'\le n$, must be carried out by deploying APCP to permute columns $k$ and $k\ma1$, then columns $k\ma1$ and $k\ma2,\dots,$ and finally columns $k'\mi1$ and $k'$. To assess the trade-offs of performing these permutations with the respective number of operations skipped in REF-ROU through the Special Case 1 adjustments, the next paragraph delves into the number of operations associated with each algorithm.

Without loss of generality, we adopt the convention that each addition, subtraction, multiplication, or division of two vector/matrix entries is considered as one operation. Permuting columns $k$ and $k\ma1$ of REF-LU$(A)$ via APCP entails $8(n\mi k)$ operations, that is, 2 multiplications, 1 addition or subtraction, and 1 division, for modifying each entry $i>k$ along each of (i) the $k$th column of $L$ (these entries entail \textit{backtracking} operations) and (ii) the $k$th row of $U$ (these entries entail \textit{row-wise switch of originating pivot} operations). Modifications made during APCP to the other matrix entries are relatively simple (e.g., sign changes) and can be assumed to take constant time (for more details on this algorithm and the two above roundoff-error-free subroutines, see \cite{esc17rou}). Conversely, calculating the $k$th column of REF-LU$(\hA)$ in REF-ROU entails $12(n\mi k)$ operations: 2 multiplications, 1 addition, and 1 division, for each entry $i>k$ along each of (i) $\zvc{k}$, (ii) the $k$th column of $\hat L$, and (iii) the diagonal of $\hat L$.

To proceed, let $\theta^1,\theta^2,\dots,\theta^p$ denote the indices of the first $p>1$ non-zeros of $\wv$, where $\theta^1<\theta^2<\dots<\theta^p$. In addition, assume that there are at least $p\mi1$ zeros among the last $n\mi p$ elements of this vector, and let $\zeta^1,\zeta^2,\dots,\zeta^{p\mi1}$ denote the indices of the first $p\mi1$ of these zeros. Based on the preceding analysis, the number of operations needed to increase $\theta_w$ from $\theta^1$ to $\theta^p$ through the requisite permutations of $\wv$ and REF-LU$(A)$ is given by\llV\lV
\begin{flalign}\label{eqn:APCP_opers}
  \sum_{j=1}^{p\mi1}\sum_{k=\theta^j}^{\zeta^{j}\mi1} 8(n-k).&
\end{flalign}
The value of the above expression must be compared with the number of operations saved by skipping columns
$\theta^1,\theta^2,\dots,\theta^{p}\mi1$ in REF-ROU, which is given by
\begin{flalign}\label{eqn:skip_opers}
  \sum_{k=\theta^1}^{\theta^{p}\mi1} 12(n-k).&
\end{flalign}
The parameter $p$ can be gradually increased and the above two expressions easily recalculated to determine precisely when the general sparsity of $\wv$ can be exploited---i.e., as long as \eqref{eqn:APCP_opers} is smaller than \eqref{eqn:skip_opers}.

Similarly, it may be possible to make efficiency gains when $\vv$ and $\wv$ are both sparse. However, while the associated permutations on REF-LU require twice the number of operations as the preceding case, the number of operations saved from skipping row and column $k$ of REF-LU$(\hA)$ increases only to $20(n-k)$ (see the proof to Theorem \ref{thm:Rank-one_Operations} for more details). It is possible that REF-ROU could be further expedited by leveraging the additional sparsity of $A$. This will be explored as part of future work.

\end{APPENDIX}
%
%
%


\bibliographystyle{informs2014}
\bibliography{References} 

\begin{thebibliography}{63}
\providecommand{\natexlab}[1]{#1}
\providecommand{\url}[1]{\texttt{#1}}
\providecommand{\urlprefix}{URL }

\bibitem[{Abbott \protect\BIBand{} Mulders(2001)}]{abb01tig}
Abbott J, Mulders T (2001) How tight is {H}adamard's bound? \emph{Experimental
  Mathematics} 10(3):331--336.

\bibitem[{Bailey \protect\BIBand{} Borwein(2015)}]{bai15hig}
Bailey DH, Borwein JM (2015) High-precision arithmetic in mathematical physics.
  \emph{Mathematics} 3(2):337--367.

\bibitem[{Bareiss(1968)}]{bar68syl}
Bareiss EH (1968) Sylvester's identity and multistep integer-preserving
  {G}aussian elimination. \emph{Mathematics of Computation} 22(103):565--578.

\bibitem[{Bareiss(1972)}]{bar72com}
Bareiss EH (1972) Computational solutions of matrix problems over an integral
  domain. \emph{IMA Journal of Applied Mathematics} 10(1):68--104.

\bibitem[{Bennett(1965)}]{ben65tri}
Bennett JM (1965) Triangular factors of modified matrices. \emph{Numerische
  Mathematik} 7(3):217--221.

\bibitem[{Choi et~al.(1990)Choi, Monma, \protect\BIBand{} Shanno}]{cho90fur}
Choi IC, Monma CL, Shanno DF (1990) Further development of a primal-dual
  interior point method. \emph{ORSA Journal on Computing} 2(4):304--311.

\bibitem[{Cook \protect\BIBand{} Steffy(2011)}]{coo11sol}
Cook W, Steffy DE (2011) Solving very sparse rational systems of equations.
  \emph{ACM Transactions on Mathematical Software (TOMS)} 37(4):39.

\bibitem[{Davis \protect\BIBand{} Hager(2001)}]{dav01mul}
Davis TA, Hager WW (2001) Multiple-rank modifications of a sparse {C}holesky
  factorization. \emph{SIAM Journal on Matrix Analysis and Applications}
  22(4):997--1013.

\bibitem[{Davis \protect\BIBand{} Hager(2005)}]{dav05row}
Davis TA, Hager WW (2005) Row modifications of a sparse {C}holesky
  factorization. \emph{SIAM Journal on Matrix Analysis and Applications}
  26(3):621--639.

\bibitem[{Deng(2010)}]{den10mul}
Deng L (2010) \emph{Multiple-rank Updates to Matrix Factorizations for
  Nonlinear Analysis and Circuit Design} (Stanford University).

\bibitem[{Diaz-Hernandez et~al.(2021)Diaz-Hernandez, Fern{\'a}ndez,
  Romano-Moreno, \protect\BIBand{} Lara}]{dia21imp}
Diaz-Hernandez G, Fern{\'a}ndez BR, Romano-Moreno E, Lara JL (2021) An improved
  model for fast and reliable harbour wave agitation assessment. \emph{Coastal
  Engineering} 104011.

\bibitem[{Drgo{\v{n}}a et~al.(2017)Drgo{\v{n}}a, Klau{\v{c}}o, Jane{\v{c}}ek,
  \protect\BIBand{} Kvasnica}]{drg17opt}
Drgo{\v{n}}a J, Klau{\v{c}}o M, Jane{\v{c}}ek F, Kvasnica M (2017) Optimal
  control of a laboratory binary distillation column via regionless explicit
  {MPC}. \emph{Computers \& Chemical Engineering} 96:139--148.

\bibitem[{Edmonds(1967)}]{edm67sys}
Edmonds J (1967) Systems of distinct representatives and linear algebra.
  \emph{Journal of Research of the National Bureau of Standards, Section B}
  71:241--245.

\bibitem[{Elble \protect\BIBand{} Sahinidis(2012)}]{elb12rev}
Elble JM, Sahinidis NV (2012) A review of the {LU} update in the simplex
  algorithm. \emph{International Journal of Mathematics in Operational
  Research} 4(4):366--399.

\bibitem[{Elsner \protect\BIBand{} Rozsa(1981)}]{els81eig}
Elsner L, Rozsa P (1981) On eigenvectors and adjoints of modified matrices.
  \emph{Linear and Multilinear Algebra} 10(3):235--247.

\bibitem[{Escobedo(2016)}]{esc16fou}
Escobedo AR (2016) \emph{Foundational Factorization Algorithms for the
  Efficient Roundoff-error-free Solution of Optimization Problems}. Ph.D.
  thesis.

\bibitem[{Escobedo \protect\BIBand{} Moreno-Centeno(2015)}]{esc15rou}
Escobedo AR, Moreno-Centeno E (2015) Roundoff-error-free algorithms for solving
  linear systems via {C}holesky and {LU} factorizations. \emph{INFORMS Journal
  on Computing} 27(4):677--689.

\bibitem[{Escobedo \protect\BIBand{} Moreno-Centeno(2017)}]{esc17rou}
Escobedo AR, Moreno-Centeno E (2017) Roundoff-error-free basis updates of {LU}
  factorizations for the efficient validation of optimality certificates.
  \emph{SIAM Journal on Matrix Analysis and Appls} 38(3):829--853.

\bibitem[{Escobedo et~al.(2018)Escobedo, Moreno-Centeno, \protect\BIBand{}
  Lourenco}]{esc18sol}
Escobedo AR, Moreno-Centeno E, Lourenco C (2018) Solution of dense linear
  systems via roundoff-error-free factorization algorithms: Theoretical
  connections and computational comparisons. \emph{ACM Transactions on
  Mathematical Software (TOMS)} 44(4):1--24.

\bibitem[{Fine \protect\BIBand{} Scheinberg(2001)}]{fin01eff}
Fine S, Scheinberg K (2001) Efficient {SVM} training using low-rank kernel
  representations. \emph{Journal of Machine Learning Research} 2(Dec):243--264.

\bibitem[{Fletcher \protect\BIBand{} Matthews(1985)}]{fle85sta}
Fletcher R, Matthews S (1985) A stable algorithm for updating triangular
  factors under a rank one change. \emph{Mathematics of Computation}
  45(172):471--485.

\bibitem[{Fletcher \protect\BIBand{} Powell(1974)}]{fle74mod}
Fletcher R, Powell MJ (1974) On the modification of ${L}{D}{L}$ factorizations.
  \emph{Mathematics of Computation} 28(128):1067--1087.

\bibitem[{Gamrath et~al.(2020)Gamrath, Anderson, Bestuzheva, Chen, Eifler,
  Gasse, Gemander, Gleixner, Gottwald, Halbig et~al.}]{scip}
Gamrath G, Anderson D, Bestuzheva K, Chen WK, Eifler L, Gasse M, Gemander P,
  Gleixner A, Gottwald L, Halbig K, et~al. (2020) The {SCIP} optimization suite
  7.0 .

\bibitem[{G{\"a}rtner \protect\BIBand{} Sch{\"o}nherr(2000)}]{gar20eff}
G{\"a}rtner B, Sch{\"o}nherr S (2000) An efficient, exact, and generic
  quadratic programming solver for geometric optimization. \emph{Proceedings of
  the Sixteenth Annual Symposium on Computational Geometry}, 110--118.

\bibitem[{Gill et~al.(1974)Gill, Golub, Murray, \protect\BIBand{}
  Saunders}]{gil74met}
Gill PE, Golub GH, Murray W, Saunders MA (1974) Methods for modifying matrix
  factorizations. \emph{Mathematics of Computation} 28(126):505--535.

\bibitem[{Gill et~al.(1987)Gill, Murray, Saunders, \protect\BIBand{}
  Wright}]{gil87mai}
Gill PE, Murray W, Saunders MA, Wright MH (1987) Maintaining {LU} factors of a
  general sparse matrix. \emph{Linear Algebra and Its Applications}
  88:239--270.

\bibitem[{Gill et~al.(1996)Gill, Saunders, \protect\BIBand{}
  Shinnerl}]{gil96sta}
Gill PE, Saunders MA, Shinnerl JR (1996) On the stability of {C}holesky
  factorization for symmetric quasidefinite systems. \emph{SIAM Journal on
  Matrix Analysis and Applications} 17(1):35--46.

\bibitem[{Gill \protect\BIBand{} Wong(2015)}]{gil15met}
Gill PE, Wong E (2015) Methods for convex and general quadratic programming.
  \emph{Mathematical Programming Computation} 7(1):71--112.

\bibitem[{Gleixner et~al.(2015)Gleixner, Miltenberger, \protect\BIBand{}
  M\"uller}]{soplex}
Gleixner A, Miltenberger M, M\"uller B (2015) So{P}lex: the sequential
  object-oriented simplex class library, version 2.2. Available at
  http://soplex.zib.de.

\bibitem[{Gleixner(2015)}]{gle15exa}
Gleixner AM (2015) \emph{Exact and Fast Algorithms for Mixed-integer Nonlinear
  Programmin}. Ph.D. thesis.

\bibitem[{Griewank \protect\BIBand{} Walther(2002)}]{gri02con}
Griewank A, Walther A (2002) On constrained optimization by adjoint based
  quasi-{N}ewton methods. \emph{Optimization Methods and Software}
  17(5):869--889.

\bibitem[{Hammarling \protect\BIBand{} Lucas(2008)}]{ham08upd}
Hammarling S, Lucas C (2008) Updating the {QR} factorization and the least
  squares problem. Technical report.

\bibitem[{Herceg et~al.(2015)Herceg, Jones, \protect\BIBand{}
  Morari}]{her15dom}
Herceg M, Jones C, Morari M (2015) Dominant speed factors of active set methods
  for fast {MPC}. \emph{Optimal Control Applications and Methods}
  36(5):608--627.

\bibitem[{Herholz \protect\BIBand{} Alexa(2018)}]{her18fac}
Herholz P, Alexa M (2018) Factor once: reusing {C}holesky factorizations on
  sub-meshes. \emph{ACM Transactions on Graphics (TOG)} 37(6):1--9.

\bibitem[{Herholz \protect\BIBand{} Sorkine-Hornung(2020)}]{her20spa}
Herholz P, Sorkine-Hornung O (2020) Sparse {C}holesky updates for interactive
  mesh parameterization. \emph{ACM Transactions on Graphics (TOG)} 39(6):1--14.

\bibitem[{Higham(2009)}]{hig09cho}
Higham NJ (2009) Cholesky factorization. \emph{Wiley Interdisciplinary Reviews:
  Computational Statistics} 1(2):251--254.

\bibitem[{Higham(2011)}]{hig11gau}
Higham NJ (2011) Gaussian elimination. \emph{Wiley Interdisciplinary Reviews:
  Computational Statistics} 3(3):230--238.

\bibitem[{Hock \protect\BIBand{} Schittkowski(1980)}]{hoc80tes}
Hock W, Schittkowski K (1980) Test examples for nonlinear programming codes.
  \emph{Journal of Optimization Theory and Applications} 30(1):127--129.

\bibitem[{Kie{\l}basi{\'n}ski \protect\BIBand{} Schwetlick(1988)}]{kie88num}
Kie{\l}basi{\'n}ski A, Schwetlick H (1988) \emph{Numerische Lineare Algebra:
  Eine Computerorientierte Einf{\"u}hrung} (Dt. Verlag d. Wiss.).

\bibitem[{Kirches et~al.(2011)Kirches, Bock, Schl{\"o}der, \protect\BIBand{}
  Sager}]{kir11fac}
Kirches C, Bock HG, Schl{\"o}der JP, Sager S (2011) A factorization with update
  procedures for a {KKT} matrix arising in direct optimal control.
  \emph{Mathematical Programming Computation} 3(4):319--348.

\bibitem[{Knuth(1981)}]{knu81art}
Knuth DE (1981) \emph{The Art of Computer Programming, Volume 2: Seminumerical
  Algorithms} (Addison-Wesley Professional), 2 edition.

\bibitem[{Lee \protect\BIBand{} Saunders(1995)}]{lee95fra}
Lee HR, Saunders BD (1995) Fraction free {G}aussian elimination for sparse
  matrices. \emph{Journal of Symbolic Computation} 19(5):393--402.

\bibitem[{Lekhovytskiy(2018)}]{lek18ada}
Lekhovytskiy DI (2018) Adaptive lattice filters for systems of space-time
  processing of non-stationary {G}aussian processes. \emph{Radioelectronics and
  Communications Systems} 61(11):477--514.

\bibitem[{Lourenco et~al.(2020)Lourenco, Chen, Moreno-Centeno,
  \protect\BIBand{} Davis}]{lou20use}
Lourenco C, Chen J, Moreno-Centeno E, Davis TA (2020) User guide for {SLIP LU},
  a sparse left-looking integer preserving {LU} factorization version 1.0.2.

\bibitem[{Lourenco et~al.(2019)Lourenco, Escobedo, Moreno-Centeno,
  \protect\BIBand{} Davis}]{lou19exa}
Lourenco C, Escobedo AR, Moreno-Centeno E, Davis TA (2019) Exact solution of
  sparse linear systems via left-looking roundoff-error-free {LU} factorization
  in time proportional to arithmetic work. \emph{SIAM Journal on Matrix
  Analysis and Appls} 40(2):609--638.

\bibitem[{Lourenco(2020)}]{lou20eff}
Lourenco CJ (2020) \emph{Efficient Algorithms for the Exact Solution of Sparse
  Linear Systems in Time Proportional to Arithmetic Work}. Ph.D. thesis.

\bibitem[{Magron et~al.(2017)Magron, Constantinides, \protect\BIBand{}
  Donaldson}]{mag17cer}
Magron V, Constantinides G, Donaldson A (2017) Certified roundoff error bounds
  using semidefinite programming. \emph{ACM Transactions on Mathematical
  Software (TOMS)} 43(4):34.

\bibitem[{Mehrotra(1992)}]{meh92def}
Mehrotra S (1992) Deferred rank one updates in o (n\_3{L}) interior point
  algorithm. \emph{Journal of the Operations Research Society of Japan}
  35(4):345--352.

\bibitem[{Minka(2013)}]{min13exp}
Minka TP (2013) Expectation propagation for approximate {B}ayesian inference.
  \emph{arXiv preprint arXiv:1301.2294} .

\bibitem[{Oh \protect\BIBand{} Hu(2018)}]{oh18mul}
Oh H, Hu Z (2018) Multiple-rank modification of symmetric eigenvalue problem.
  \emph{MethodsX} 5:103--117.

\bibitem[{Ojeda et~al.(2008)Ojeda, Suykens, \protect\BIBand{}
  De~Moor}]{oje08low}
Ojeda F, Suykens JA, De~Moor B (2008) Low rank updated {LS-SVM} classifiers for
  fast variable selection. \emph{Neural Networks} 21(2-3):437--449.

\bibitem[{Olszanskyj et~al.(1994)Olszanskyj, Lebak, \protect\BIBand{}
  Bojanczyk}]{ols94ran}
Olszanskyj SJ, Lebak JM, Bojanczyk AW (1994) Rank-k modification methods for
  recursive least squares problems. \emph{Numerical Algorithms} 7(2):325--354.

\bibitem[{Pan(2015)}]{pan15var}
Pan PQ (2015) The variant of the face algorithm is unstable. \emph{preprint} .

\bibitem[{Pan(2020)}]{pan20new}
Pan PQ (2020) A new face algorithm using {LU} factorization for linear
  programming. \emph{preprint} .

\bibitem[{Puranik \protect\BIBand{} Sahinidis(2017)}]{pur17bou}
Puranik Y, Sahinidis NV (2017) Bounds tightening based on optimality conditions
  for nonconvex box-constrained optimization. \emph{Journal of Global
  Optimization} 67(1-2):59--77.

\bibitem[{Sarra(2011)}]{sar11rad}
Sarra SA (2011) Radial basis function approximation methods with extended
  precision floating point arithmetic. \emph{Engineering Analysis with Boundary
  Elements} 35(1):68--76.

\bibitem[{Sch{\"o}nhage \protect\BIBand{} Strassen(1971)}]{sch71sch}
Sch{\"o}nhage DDA, Strassen V (1971) Schnelle multiplikation grosser zahlen.
  \emph{Computing} 7(3-4):281--292.

\bibitem[{Seeger(2008)}]{see08low}
Seeger M (2008) Low rank updates for the {C}holesky decomposition. Technical
  report.

\bibitem[{Seeger et~al.(2007)Seeger, Steinke, \protect\BIBand{}
  Tsuda}]{see07bay}
Seeger M, Steinke F, Tsuda K (2007) Bayesian inference and optimal design in
  the sparse linear model. \emph{Artificial Intelligence and Statistics},
  444--451 (PMLR).

\bibitem[{Stange et~al.(2007)Stange, Griewank, \protect\BIBand{}
  Bollh\"{o}fer}]{sta07eff}
Stange P, Griewank A, Bollh\"{o}fer M (2007) On the efficient update of
  rectangular {LU}-factorizations subject to low rank modifications.
  \emph{Electronic Transactions on Numerical Analysis} 26:161--177.

\bibitem[{Weber et~al.(2019)Weber, Sager, \protect\BIBand{}
  Gleixner}]{web19sol}
Weber T, Sager S, Gleixner A (2019) Solving quadratic programs to high
  precision using scaled iterative refinement. \emph{Mathematical Programming
  Computation} 11(3):421--455.

\bibitem[{Wright et~al.(1999)Wright, Nocedal et~al.}]{wri99num}
Wright S, Nocedal J, et~al. (1999) Numerical optimization. \emph{Springer
  Science} 35(67-68):7.

\bibitem[{Wunderling(1996)}]{wun96par}
Wunderling R (1996) \emph{Paralleler und Objektorientierter
  {S}implex-{A}lgorithmus}. Ph.D. thesis, Technische Universit{\"a}t Berlin,
  \url{http://www.zib.de/Publications/abstracts/TR-96-09/}.

\end{thebibliography}

\end{document}